\documentclass[11pt,a4paper]{article}

\RequirePackage{amsmath}%
\RequirePackage{amsthm}%
\usepackage{amsfonts}%
\usepackage{amssymb}%
\usepackage{graphicx}%
\usepackage{cite}

\newtheorem{theorem}{Theorem}
\newtheorem{lemma}{Lemma}

\newcommand{ \Null}{{\rm Null}}

\renewcommand{\Re}{\mathop{\mathrm{Re}}}
\renewcommand{\Im}{\mathop{\mathrm{Im}}}

\renewcommand{\i}{\mathrm{i}}
\newcommand{\Res}{\mathop{\mathrm{Res}}}

\newcommand{\bI}{{\bf I}}

\newcommand{\bM}{{\bf M}}
\newcommand{\bN}{{\bf N}}
\newcommand{\bJ}{{\bf J}}

\begin{document}

\title{Numerical conformal mapping via a boundary integral equation with the adjoint generalized Neumann kernel}

\author{Mohamed M.S. Nasser$^{\rm a,b}$, Ali H.M. Murid$^{\rm c}$ and Ali W.K. Sangawi$^{\rm d}$  \\[10pt] 
$^{\rm a}$Department of Mathematics, Faculty of Science, \\
King Khalid University, P. O. Box 9004, Abha, Saudi Arabia.\\
$^{\rm b}$Department of Mathematics, Faculty of Science, \\
Ibb University, P.O.Box 70270, Ibb, Yemen.\\
$^{\rm c}$Department of Mathematical Sciences, Faculty of Science,\\ 
Universiti Teknologi Malaysia, 81310 UTM Johor Bahru, Johor, Malaysia.\\
$^{\rm d}$Department of Mathematics, College of Science, \\
University of Sulaimani, 46001 UOS, Sulaimani, Kurdistan, Iraq.\\
{\tt mms\_nasser@hotmail.com, alihassan@utm.my, alisangawi2000@yahoo.com}\\
}

\date{}
\maketitle


\begin{center}
\begin{quotation}
{\noindent {\bf Abstracts.\;\;}%
This paper presents a new uniquely solvable boundary integral equation for computing the conformal mapping, its derivative and its inverse from bounded multiply connected regions onto the five classical canonical slit regions. The integral equation is derived by reformulating the conformal mapping as an adjoint Riemann-Hilbert problem. From the adjoint Riemann-Hilbert problem, we derive a boundary integral equation with the adjoint generalized Neumann kernel for the derivative of the boundary correspondence function $\theta'$. Only the right-hand side of the integral equation is different from a canonical region to another. The function $\theta'$ is integrated to obtain the boundary correspondence function $\theta$. The integration constants as well as the parameters of the canonical region are computed using the same uniquely solvable integral equation. 

A numerical example is presented to illustrate the accuracy of the proposed method. 
}%

\end{quotation}
\end{center}

\begin{center}
\begin{quotation}
{\noindent {\bf Keywords.\;\;}%
Numerical conformal mapping; Multiply connected regions; Generalized Neumann
kernel; Riemann-Hilbert problem.}%
\end{quotation}
\end{center}

\begin{center}
\begin{quotation}
{\noindent {\bf MSC.\;\;} 30C30; 30E25; 65E05. }
\end{quotation}
\end{center}

\section{Introduction}
\label{sc:int}

The classical canonical slit regions for multiply connected regions are: an annulus with concentric circular slit region, a disk with concentric circular slit region, the circular slit region, the radial slit region, and the parallel slit region (see~\cite{ber,neh}). Assume that $\Omega$ is any one of these five canonical regions. It is well-known that a conformal mapping $w=\omega(z)$ form any bounded multiply connected region $G$ in the $z$-plane onto the canonical slit region $\Omega$ in the $w$-plane always exist. Several numerical methods have been proposed for computing such mapping function~\cite{ama,ama2,cro,del,kok,may,nas-fun,nas-siam,nas-jmaa,okab,rei,san-dc,san-pr,san-an,san-cr,san-rd,san-ans,san-sps,yun,yun-bul,yun-prsa}. Most of these numerical methods can be used to calculate the mapping function onto only one canonical region~\cite{may,san-dc,san-pr,san-an,san-cr,san-rd,san-ans,san-sps,yun}. Other methods can be used to compute the mapping function, in a unified way, onto two canonical regions~\cite{ama,del,kok,okab,rei} or three canonical regions~\cite{ama2}. Only the method presented in~\cite{nas-fun,nas-siam,nas-jmaa} can be used to compute the mapping function onto the five canonical regions in a unified way. The method is based on a boundary integral equation with the generalized Neumann kernel

Most of the above numerical methods require solving Riemann-Hilbert (RH) problems~\cite{ama,ama2,del,may,nas-fun,nas-siam,nas-jmaa,nas-jmaa13,okab}. The iterative method in~\cite{wegf} requires at each iterative step the solution of a RH problem. In~\cite{ama,ama2,may,okab}, the mapping function is expressed in terms of a solution of a modified Dirichlet problem which is a special case of the RH problem. 

In~\cite{nas-fun,nas-siam}, the mapping function $\omega(z)$ onto the above five classical canonical slit regions was reformulated in terms of an auxiliary function $f(z)$ which is a solution of a RH problem. Then, by the properties of the mapping function on the boundary $\Gamma:=\partial G$, the function $f(z)$ is written as a solution of a RH problem with the coefficient function $A(t)=\eta(t)$, where $\eta(t)$ is a parametrization of the boundary $\Gamma$. The RH problem is solved using a uniquely solvable boundary integral equation with the generalized Neumann kernel. By solving the integral equation, we obtain the boundary values of the auxiliary function $f(z)$ and the parameters of the canonical regions, and hence we obtain the boundary values of the mapping function $\omega(z)$. The values of the mapping function $\omega(z)$ for interior points are computed by the Cauchy integral formula. The method can be used for regions $G$ with smooth or piecewise smooth boundaries. 

A formulation of the conformal mapping from unbounded multiply connected circular regions onto the radial or circular slit domain, or to domains with both radial and circular slit in terms of a RH problem is given in~\cite{del}. The same approach was used in~\cite{mit} for conformal mapping from multiply connected circular regions onto the circular slit domain. These methods can be used only if the boundaries of the original region $G$ are circles.  The function $\omega'(z)/\omega(z)$ was reformulated in terms of an auxiliary function which, by properties of the function $\omega'(z)/\omega(z)$ on the boundary $\Gamma$, is a solution of a RH problem. Since the boundaries of the original region $G$ are circles, the coefficient function of the RH problem formulated in~\cite{del,mit} is equivalent to the function $\tilde A(t)=\dot\eta(t)$, where $\tilde A$ is the adjoint of the function $A(t)=1$ (see~\cite[Eq.~(11)]{wegm} for the definition of the adjoint function). Hence, the RH problem formulated in~\cite{del,mit} is the adjoint of the RH problem with the coefficient function $A(t)=1$.

This paper presents a combination of the approaches presented in~\cite{del,mit} and \cite{nas-fun,nas-siam,nas-jmaa,nas-jmaa13}. In formulating the mapping function as a RH problem, we shall use the approach presented~\cite{del,mit}. For solving the RH problem, we shall use the boundary integral equation approach as in~\cite{nas-fun,nas-siam,nas-jmaa,nas-jmaa13}. 

Firstly, we used the approach presented in~\cite{del,mit}, with slight modifications, to reformulate the mapping function from bounded multiply connected regions onto the above five classical canonical slit regions in terms of a RH problems. The mapping function will be written in terms of an auxiliary function $F(z)$. We define also an auxiliary function $f(z)$ in terms of the function $F(z)$ and its derivative $F'(z)$. For the first four canonical regions, the properties of the function $\omega'(z)/\omega(z)$ on the boundary $\Gamma$ implies that the function $f(z)$ is a solution of a RH problem with the coefficient function $\tilde A(t)=\dot\eta(t)/\eta(t)$. For the last canonical region, we get the same RH problem form the properties of the function $\omega'(z)$ on the boundary $\Gamma$. The RH problem which will be formulated in this paper is the adjoint of the RH problem formulated in~\cite{nas-fun,nas-siam}. Only the right-hand side of the RH problem is different from a canonical region to another. 
 
Secondly, based on the results presented in~\cite{wegm}, we shall use the formulated RH problem to derive a boundary integral equation for the derivative $\theta'$ of the boundary correspondence function. The kernel of the derived integral equation is the adjoint of the generalized Neumann kernel in~\cite{nas-fun,nas-siam}. The derived integral equation has appeared in~\cite{nas-lap} for solving the Dirichlet problem and the Neumann problem in multiply connected regions. It is not uniquely solvable. However, based on the properties of the function $\theta'$ and in view of the results presented in~\cite{nas-lap}, the integral equation can be modified to obtain a uniquely solvable boundary integral equation for computing the function $\theta'$. Only the right-hand side of the integral equation is different from a canonical region to another. The boundary correspondence function $\theta$ will be computed as an anti-derivative of the function $\theta'$. The integration constants as well as the parameters of the canonical regions will be calculated by solving the same uniquely solvable integral equation. Hence, we obtain the boundary values of the mapping function $\omega(z)$ and the boundary values of its derivative $\omega'(z)$. The values of the mapping function $w=\omega(z)$ and its derivative $\omega'(z)$ for $z\in G$ will be calculated by means of the Cauchy integral formula. The presented method can also be used to compute the inverse mapping function $z=\omega^{-1}(w)$ from $\Omega$ onto $G$. The values of the inverse mapping function $z=\omega^{-1}(w)$ will be calculated from the boundary correspondence function $\theta$ and its derivative $\theta'$ using the same approach used in~\cite[p.~380]{hen} for simply connected regions. 

Other integral equations for numerical computing of conformal mapping of multiply connected regions have been derived in~\cite{san-dc,san-pr,san-an,san-cr,san-rd,yun}. The approach used in~\cite{san-dc,san-pr,san-an,san-cr,san-rd,yun} for deriving the integral equations is completely different from the approach used in this paper or in~\cite{nas-fun,nas-siam,nas-jmaa,nas-jmaa13}. It is an extension of the approach used in~\cite{mur,raz} to derive boundary integral equation for conformal mapping of simply connected regions. It has been used also in~\cite{mur-hu} to derive an integral equation for conformal mapping of bounded multiply connected regions onto an annulus with circular slits.

The method presented in~\cite{san-dc,san-pr,san-an,san-cr,san-rd,yun} depends on three integral equations. The first two boundary integral equations were derived from certain boundary relationships on the boundary $\Gamma$. The third integral equation is the boundary integral equation with adjoint generalized Neumann kernel which has been derived in~\cite{nas-lap}. The first integral equation is used to calculate the derivative $\omega'(z)$. However, the kernel of this integral equation contains the parameters of the canonical region which should be calculated first using the third integral equation. The second integral equation is used to calculate $\theta'$. The kernel of the second integral equation is the adjoint Neumann kernel which is different from the adjoint generalized Neumann kernel. Then, the boundary values of the mapping function $\omega(z)$ are computed from a certain boundary relationship that relates the mapping function $\omega(z)$ on the boundary $\Gamma$ with its derivative $\omega'(z)$ , $\theta'$, and the parameters of the canonical regions.

\section{Notations and auxiliary material}
\label{sc:aux}

\begin{figure} %
\centerline{\scalebox{0.6}[0.30]{\includegraphics{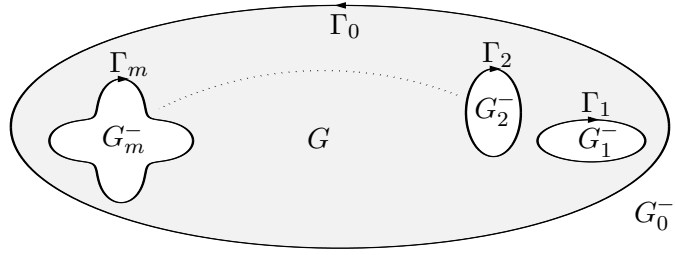}}}
 \vskip-3.35cm \noindent\hspace{7.2cm} $\Gamma_0$
 \vskip-0.10cm   \noindent\hspace{9.2cm} $\Gamma_2$
 \vskip-0.35cm  \noindent\hspace{4.3cm}  $\Gamma_m$
 \vskip+0.05cm  \noindent\hspace{10.5cm}  $\Gamma_1$
 \vskip-.40cm   \noindent\hspace{9.1cm} $G_2^-$
 \vskip-.10cm  \noindent\hspace{4.2cm}  $G^-_m$
 \vskip-0.45cm  \noindent\hspace{6.9cm} $G$ \hspace{3.0cm} $G_1^-$
 \vskip+0.5cm  \noindent\hspace{11.2cm} $G_0^-$
 \vskip+0.25cm
\caption{\rm The bounded multiply connected region $G$.} \label{f:bd}
\end{figure}

We consider bounded multiply connected regions $G$ of connectivity $m+1\ge1$ with boundary $\Gamma=\cup_{j=0}^{m}\Gamma_j$ consisting of $m+1$ smooth closed Jordan curves $\Gamma_j$, $j=0,1,2,\ldots,m$. The outer boundary curve $\Gamma_0$ is oriented counterclockwise and the inner boundaries $\Gamma_1,\dots,\Gamma_m$ are oriented clockwise. The complement $G^-:=\overline{\mathbb{C}} \setminus \overline{G}$ consists of $m$ bounded simply connected components $G_j$ interior to $\Gamma_j$, $j=1,2, \ldots, m$, and an unbounded simply connected component $G_0$ exterior to $\Gamma_0$ (see Figure~\ref{f:bd}). The curve $\Gamma_j$ is parametrized by a $2\pi$-periodic twice continuously differentiable complex function $\eta_j(t)$ with
non-vanishing first derivative
\begin{equation}\label{e:par}
\dot\eta_j(t)=d\eta_j(t)/dt\ne 0, \quad t\in J_j:=[0,2\pi], \quad j=0,1,\ldots,m.
\end{equation}
The total parameter region $J$ is the disjoint union of the
intervals $J_j$. We define a parametrization of the
whole boundary $\Gamma$ as the complex function $\eta$ defined on $J$ by
\begin{equation}\label{e:eta}
\eta(t):= \left\{
\begin{array}{l@{\hspace{1cm}}l}
\eta_0(t), &  t\in J_0,	\\
\vdots     &              \\
\eta_m(t), &  t\in J_m. 	\\
\end{array}\right.
\end{equation}

Let $H$ be the space of all real H\"older continuous $2\pi-$periodic functions $\phi(t)$ of the parameter $t$ on $J_j$ for $j=0,\ldots, m$, i.e.,
\[
\phi(t) = \left\{
\begin{array}{l@{\hspace{0.5cm}}l}
 \phi_0(t),     & t\in J_0, \\
  \vdots       & \\
 \phi_m(t),     & t\in J_m, \\
\end{array}%
\right.
\]
with real H\"older continuous $2\pi-$periodic functions $\phi_0,\ldots,\phi_m$. In view of the smoothness of $\eta$, a real H\"older continuous function $\hat\phi$ on $\Gamma$ can be interpreted via $\phi(t):=\hat\phi(\eta(t))$, $t\in J$, as a function $\phi\in H$; and vice versa. Let $S$ be the subspace of $H$ that consists of real piecewise constant functions of the form
\[
h(t) = \left\{
\begin{array}{l@{\hspace{0.5cm}}l}
 h_0,     & t\in J_0, \\
  \vdots       & \\
 h_m,     & t\in J_m, \\
\end{array}%
\right.
\]
with real constants $h_0,\ldots,h_m$. For simplicity, the piecewise constant function $h$ will be denoted by
\[
h(t)=(h_0,\ldots,h_m).
\]

Let $w=\omega(z)$ be the mapping function from the region $G$ onto $\Omega$ where $\Omega$ is any of the above canonical regions and let $L:=\partial\Omega$. Let also $L$ be parametrized by 
\begin{equation}\label{e:zeta}
\zeta(t), \qquad t\in\hat J,
\end{equation}
where $\hat J$ is a suitable parameter region. The function $\zeta(s)$ contains unknown parameters known as the parameters of the canonical region. The mapping function $w=\omega(z)$ is completely described by its boundary values which can be described by
\begin{equation}\label{e:theta}
\omega(\eta(t))=\zeta(\theta(t)), \quad t\in J,
\end{equation}
where the function $\theta(t)$ is called the \emph{boundary correspondence function}. It follows from~(\ref{e:theta}) that determining the function $\theta(t)$ and the parameters of the canonical region completely determine the boundary values of the mapping function $\omega(z)$. 

By differentiation both sides of~(\ref{e:theta}) with respect to the parameter $t$, we obtain
\begin{equation}\label{e:theta'}
\dot\eta(t)\omega'(\eta(t))=\theta'(t)\zeta(\theta(t)), \quad t\in J.
\end{equation}
Hence, the determining the function $\theta(t)$ and its derivative $\theta'(t)$ determine the boundary values of the derivative $\omega'(z)$ of the mapping function.

\section{Generalized Neumann kernel}
\label{sc:gnk}

Let $A$ be the continuously differentiable complex function defined by
\begin{equation}\label{e:A}
A(t) := \eta(t), \quad t\in J.
\end{equation}
We define two real kernels
\begin{eqnarray}
\label{e:N}
 N(s,t) &:=&  \frac{1}{\pi}\Im\left(
 \frac{A(s)}{A(t)}\frac{\dot\eta(t)}{\eta(t)-\eta(s)}\right), \\
\label{e:M}
 M(s,t) &:=&  \frac{1}{\pi}\Re\left(
 \frac{A(s)}{A(t)}\frac{\dot\eta(t)}{\eta(t)-\eta(s)}\right).
\end{eqnarray}
The kernel $N(s,t)$ is known as the \emph{generalized Neumann kernel} formed with $A$. The kernel $N$ is continuous and the kernel $M$ has a cotangent singularity type (see~\cite{wegm} for more details). Thus, the operator 
\begin{equation}\label{e:opN}
  \bN \mu(s) := \int_J N(s,t) \mu(t) dt, \quad s\in J,
\end{equation}
is a Fredholm integral operator and the operator 
\begin{equation}\label{e:opM}
  \bM\mu(s) := \int_J  M(s,t) \mu(t) dt, \quad s\in J,
\end{equation}
is a singular integral operator.

The function $\tilde A$ defined by
\begin{equation}\label{e:tA}
\tilde A(t) =\frac{\dot\eta(t)}{A(t)} =\frac{\dot\eta(t)}{\eta(t)}
\end{equation}
is known as the \emph{adjoint function} to the function $A$. Then, the generalized Neumann kernel $\tilde N$ formed with $\tilde A$ is defined by
\begin{equation}\label{e:tN}
 \tilde N(s,t) :=  \frac{1}{\pi}\Im\left(
 \frac{\tilde A(s)}{\tilde A(t)}\frac{\dot\eta(t)}{\eta(t)-\eta(s)}\right).
\end{equation}
We define also the real kernel $\tilde M$ by
\begin{equation}\label{e:tM}
 \tilde M(s,t) :=  \frac{1}{\pi}\Re\left(
 \frac{\tilde A(s)}{\tilde A(t)} \frac{\dot\eta(t)}{\eta(t)-\eta(s)}\right).
\end{equation}
Then,
\begin{equation}\label{e:N*=-tN}
\tilde N(s,t)=-N^\ast(s,t)\quad {\rm and}\quad \tilde M(s,t)=-M^\ast(s,t)
\end{equation}
where $N^\ast(s,t)=N(t,s)$ is the adjoint kernel of the generalized Neumann kernel $N(s,t)$ and $M^\ast(s,t)=M(t,s)$ is the adjoint kernel of the kernel $M(s,t)$  (see~\cite{nas-lap,wegm}). Let the Fredholm operator $\tilde\bN$ and the singular operator $\tilde\bM$ be defined as in~(\ref{e:opN}) and~(\ref{e:opM}). Then~(\ref{e:N*=-tN}) implies that
\begin{equation}\label{e:bN-bN*}
\bN^\ast=-\tilde\bN, \quad \bM^\ast=-\tilde\bM,
\end{equation}
where $\bN^\ast$ and $\bM^\ast$ are the adjoint operators to the operators $\bN$ and $\bM$ respectively.

Finally, we define an integral operator $\bJ$ by
\begin{equation}\label{e:bJ}
\bJ\mu(s):=\int_J \frac{1}{2\pi}\sum_{j=0}^m \chi^{[j]}(s)\chi^{[j]}(t)\mu(t)dt,
\end{equation}
where $\chi^{[j]}$ is the piecewise constant functions defined on $J$ by
\begin{equation}\label{e:chi}
\chi^{[j]}(t):= \left\{
\begin{array}{l@{\hspace{1cm}}l}
1, &  \mbox{if }\,t\in    J_j,   \\
0, &  \mbox{if }\,t\notin J_j,   \\
\end{array}\right.
\end{equation}
$j=0,1,\ldots,m$. Then, by~\cite[Theorem~4]{nas-lap}, we have
\begin{equation}\label{e:null-N*+J}
\Null(\bI+\bN^\ast+\bJ)=\{0\}.
\end{equation}
For $j=0,\ldots,m$, let $\phi^{[j]}$ be the unique solution of the integral equation
\begin{equation}\label{e:ie-phi-j}
(\bI+\bN^\ast+\bJ)\phi^{[j]}=-\chi^{[j]}.
\end{equation}
Thus, we have from~\cite{nas-lap} the following theorem.

\begin{theorem}\label{T:ie-hp}
Let $\gamma,\mu\in H$ and $h,p\in S$ such that 
\begin{equation}\label{e:Af}
Af=\gamma+h+\i[\mu+\nu]
\end{equation}
are boundary values of an analytic function $f(z)$ in $G$. Then the functions $h=(h_0,\ldots,h_m)$ and $\nu=(\nu_0,\ldots,\nu_m)$ are given by
\begin{eqnarray}
\label{e:h=}
h&=&\sum_{k=0}^{m}\left(\gamma,\phi^{[k]}\right)\chi^{[k]}, \\
\label{e:p=}
\nu&=&\sum_{k=0}^{m}\left(\mu,\phi^{[k]}\right)\chi^{[k]}. 
\end{eqnarray}
\end{theorem}
\noindent%
{\bf Proof.} \\
The formula~(\ref{e:h=}) follows from~\cite[Theorem~5]{nas-lap}. For the formula~(\ref{e:p=}), the function 
\[
\hat f(z):=-\i f(z)
\]
is analytic function in $G$ and has the boundary values
\begin{equation}\label{e:Af-hf}
A\hat f=\mu+\nu+\i(-\gamma-h).
\end{equation}
Then the formula~(\ref{e:p=}) follows from~\cite[Theorem~5]{nas-lap}.
\hfill $\Box$

\begin{theorem}\label{T:ie}
Let $\upsilon,\varphi,\psi,\phi\in H$, $f(z)$ be analytic function in $G$ and $g(z)$ be analytic function in $G^-$ with $g(\infty)=0$ such that the boundary values of the functions $f$ and $g$ are given by
\begin{equation}\label{e:tAf+tAg}
\tilde A(t) f(\eta(t))+\tilde A(t) g(\eta(t))=\upsilon+\i\varphi,
\end{equation}
where 
\begin{equation}\label{e:Jmu=th}
\bJ\varphi=\tilde h=(\tilde h_0,\ldots,\tilde h_m)
\end{equation}
is given function. Let also the boundary values of the function $g$ are given by
\begin{equation}\label{e:tAg}
\tilde A(t) g(\eta(t))=\psi+\i\phi.
\end{equation}
Then the function $\varphi$ is the unique solution of the integral equation
\begin{equation}\label{e:ie}
(\bI+\bN^\ast+\bJ)\varphi=\bM^\ast\upsilon+2\phi+\tilde h.
\end{equation}
\end{theorem}
\noindent%
{\bf Proof.} \\
It follows from~(\ref{e:tAf+tAg}) and from~(\ref{e:tAg}) that the boundary values of the function $f$ are given by
\begin{equation}\label{e:tAf}
\tilde A(t) f(\eta(t))=(\upsilon(t)-\psi(t))+\i(\varphi(t)-\phi(t)),
\end{equation}
which implies that the function $f(z)$ is a solution of the interior adjoint Riemann-Hilbert problem
\begin{equation}\label{e:rhp-i}
\Re[\tilde A(t) f(\eta(t))]=\upsilon(t)-\psi(t).
\end{equation}
Then, in view of~(\ref{e:N*=-tN}) and~(\ref{e:tAf}), it follows from~\cite[Theorem~2(c)]{wegm} that the function $\varphi-\phi$ satisfies the integral equation
\begin{equation}\label{e:ie-tAf}
(\bI+\bN^\ast)(\varphi-\phi)=\bM^\ast(\upsilon-\psi).
\end{equation}
Similarly, it follows from~(\ref{e:tAg}) that the function $g(z)$ is a solution of the exterior adjoint Riemann-Hilbert problem
\begin{equation}\label{e:rhp-e}
\Re[\tilde A(t) g(\eta(t))]=\psi(t).
\end{equation}
Then, in view of~(\ref{e:N*=-tN}) and~(\ref{e:tAg}), it follows from~\cite[Theorem~2(d)]{wegm} that the function $\phi$ satisfies the integral equation
\begin{equation}\label{e:ie-tAg}
(\bI-\bN^\ast)\phi=-\bM^\ast\psi.
\end{equation}
Subtracting~(\ref{e:ie-tAg}) form~(\ref{e:ie-tAf}) yields the integral equation 
\begin{equation}\label{e:ie-2}
(\bI+\bN^\ast)\varphi=2\phi+\bM^\ast\upsilon.
\end{equation}
By adding~(\ref{e:Jmu=th}) to~(\ref{e:ie-2}), we obtain~(\ref{e:ie}). 
\hfill $\Box$ 

The Riemann-Hilbert problem~(\ref{e:rhp-i}) is the adjoint of the Riemann-Hilbert problem~(11) in~\cite{nas-siam}. The operator $\bN^\ast$ is also the adjoint of the operator $\bN$ in~\cite{nas-siam}. In this paper, the integral equation~(\ref{e:ie}) will be used to derive an integral equation for $\theta'$. Although, the operator $\bM^\ast$ appears on the right-hand side of the integral equation~(\ref{e:ie}), the term $\bM^\ast\upsilon$ will vanish for the integral equation for $\theta'$ as we shall see in the forthcoming sections.

\begin{lemma}\label{L:cau}
If $f(z)$ is an analytic in $G$ with a simple pole at $z=0$ with
\begin{equation}\label{e:Resf}
c=\Res_{z=0} \omega(z),
\end{equation}
then the functions $f(z)$ and $f'(z)$ can be computed for $z\in G$ by
\begin{eqnarray}
\label{e:cau-b1}
f(z)&=&\frac{c}{z} +\frac{1}{2\pi\i} \int_{\Gamma} \frac{f(\eta)}{\eta-z}d\eta, \\
\label{e:cau-b1'}
f'(z)&=&-\frac{c}{z^2} +\frac{1}{2\pi\i} \int_{\Gamma} \frac{f'(\eta)}{\eta-z}d\eta. 
\end{eqnarray}
\end{lemma} 
\noindent%
{\bf Proof.} 
The function 
\begin{equation}\label{e:hf-f}
\hat f(z):=f(z)-\frac{c}{z}
\end{equation}
and its derivative
\begin{equation}\label{e:hf-f'}
\hat f'(z):=f'(z)+\frac{c}{z^2}
\end{equation}
are analytic in $G$. Hence, by the Cauchy integral formula, we have
\begin{equation}\label{e:hf-f2}
\hat f(z) =\frac{1}{2\pi\i} \int_{\Gamma} \frac{\hat f(\eta)}{\eta-z}d\eta = \frac{1}{2\pi\i} \int_{\Gamma} \frac{f(\eta)}{\eta-z}d\eta -\frac{1}{2\pi\i} \int_{\Gamma} \frac{1}{\eta-z} \frac{c}{\eta}d\eta.
\end{equation}
and
\begin{equation}\label{e:hf-f'2}
\hat f'(z) =\frac{1}{2\pi\i} \int_{\Gamma} \frac{\hat f'(\eta)}{\eta-z}d\eta = \frac{1}{2\pi\i} \int_{\Gamma} \frac{f'(\eta)}{\eta-z}d\eta +\frac{1}{2\pi\i} \int_{\Gamma} \frac{1}{\eta-z} \frac{c}{\eta^2}d\eta.
\end{equation}
Since
\[
\frac{1}{2\pi\i} \int_{\Gamma} \frac{1}{\eta-z} \frac{1}{\eta}d\eta=0\quad \text{and}\quad
\frac{1}{2\pi\i} \int_{\Gamma} \frac{1}{\eta-z} \frac{1}{\eta^2}d\eta=0,
\]
then it follows from~(\ref{e:hf-f}), (\ref{e:hf-f'}), (\ref{e:hf-f2}), and~(\ref{e:hf-f'2}) that the functions $f(z)$ and $f'(z)$ are given by~(\ref{e:cau-b1}) and~(\ref{e:cau-b1'}).
\hfill $\Box$

\section{Computing the function $\theta$}
\label{sc:thet}

Assume that $\theta(t)$, $t\in J$, is the boundary correspondence function of the mapping function from the region $G$ onto any of the canonical slit regions listed above. In this paper, we shall calculate the derivative $\theta'(t)$ using a boundary integral equation with the adjoint generalized Neumann kernel. More precisely, we shall show that the function $\theta'$ is the unique solution of the integral equation
\begin{equation}\label{e:thet'-ie}
(\bI+\bN^\ast+\bJ)\theta'=\hat\phi,
\end{equation}
where only the the right-hand side term $\hat\phi$ is different from a canonical region to another. 

For $k=0,1,\ldots,m$, the function $\theta_k(t)$ can be calculated from $\theta'_k(t)$ by 
\begin{equation}\label{e:thet-rho}
\theta_k(t)=\int\theta'_k(t)dt+c_k=:\rho_k(t)+c_k, \quad t\in J_k,
\end{equation}
where $c_k$ is undetermined real constant and the real function $\rho_k(t)$ is defined by
\begin{equation}\label{e:rho-k}
\rho_k(t):=\int\theta'_k(t)dt, \quad t\in J_k.
\end{equation}
The undetermined constants $c_k$ in~(\ref{e:thet-rho}) will be calculated using the method explained in Theorem~\ref{T:ie-hp} and the function $\rho_k(t)$ will be calculated using Fourier series. The function $\theta_k(t)$ is not necessary a $2\pi$-periodic. However, the derivative function $\theta'_k(t)$ is $2\pi$-periodic. Thus, the function $\theta'_k(t)$ can be represented by a Fourier series 
\begin{equation}\label{e:thet'-F}
\theta'_k(t)=a_0^{[k]}+\sum_{j=1}^{\infty}a_j^{[k]}\cos jt +\sum_{j=1}^{\infty}b_j^{[k]}\sin jt, \quad t\in J_k.
\end{equation}
Hence the function $\rho_k(t)$ has the Fourier series representation
\begin{equation}\label{e:rho-F}
\rho_k(t) =a_0^{[k]}t+\sum_{j=1}^{\infty}\frac{a_j^{[k]}}{j}\sin jt -\sum_{j=1}^{\infty}\frac{b_j^{[k]}}{j}\cos jt, \quad t\in J_k.
\end{equation}

\section{An annulus with circular slit region}

This canonical region $\Omega$ consists of a circular ring centered at the origin slit along $m-1$ arcs of circles. We assume that $\omega$ maps the curve $\Gamma_0$ onto the unit circle $|w|=1$, the curve $\Gamma_1$ onto the circle $|w|=R_1$ and the curves $\Gamma_j$, $j=2,3,\ldots,m$, onto slits on the circles $|w|=R_j$, where $R_1,\ldots,R_m$ are undetermined real constants. Then, the boundary values of the mapping function $\omega$ are given by
\begin{equation}\label{e:rn-w}
\omega(\eta(t))=R(t)e^{\i\theta(t)}
\end{equation}
where $\theta(t)$ is the boundary correspondence function and $R(t)=(1,R_1,\ldots,R_m)$. Thus, by differentiating both sides of~(\ref{e:rn-w}), we obtain
\begin{equation}\label{e:rn-w'}
\dot\eta(t)\frac{\omega'(\eta(t))}{\omega(\eta(t))}=\i\theta'(t).
\end{equation}

The function $\omega(z)$ is uniquely determined by assuming
\begin{equation}\label{e:ann-b-cd1}
\omega(0)>0.
\end{equation}
Thus the function $\omega(z)$ can be expressed in the form
\begin{equation}\label{e:rn-wb}
\omega(z)= c\,\left(1-\frac{z}{z_1}\right)\,
e^{zF(z)}
\end{equation}
where $c:=\omega(0)$ is an undetermined positive real constant, $z_1$ is a fixed point in $G_1$ and $F(z)$ is an auxiliary unknown function. Hence
\begin{equation}\label{e:rn-wb1}
\eta(t)F(\eta(t))+\ln c+\log\left(1-\frac{\eta(t)}{z_1}\right) =\log(\omega(\eta(t))).
\end{equation} 
Thus, by differentiating both sides of~(\ref{e:rn-wb1}) and using~(\ref{e:rn-w'}), we obtain
\begin{equation}\label{e:rn-b-bd}
\dot\eta(t)\left(F(\eta(t))+\eta(t)F'(\eta(t))\right) +\dot\eta(t)\frac{1}{\eta(t)-z_1} =\i\theta'(t).
\end{equation}

The function $f(z)$ defined in $G$ by
\begin{equation}\label{e:rn-f}
f(z):= zF(z)+z^2F'(z)+\frac{z}{z-z_1}, 
\end{equation}
and the function $g(z)$ defined in $G^-$ by $g(z):=0$ satisfy the assumptions of Theorem~\ref{T:ie} and their boundary values satisfy~(\ref{e:tAf+tAg}) with 
\begin{equation}\label{e:rn-b-gam-mu}
\upsilon(t)=0\quad{\rm and}\quad\varphi(t)=\theta'(t).
\end{equation}
Since the image of the curve $\Gamma_0$ is counterclockwise oriented, the image of the curve $\Gamma_1$ is clockwise oriented and the images of the curves $\Gamma_j$, $j=2,\ldots,m$, are slits so we have $\theta_0(2\pi)-\theta_0(0)=2\pi$, $\theta_1(2\pi)-\theta_1(0)=-2\pi$ and $\theta_j(2\pi)-\theta_j(0)=0$. Hence the function $\tilde h(t)$ in~(\ref{e:ie}) is given by
\begin{equation}\label{e:rn-b-th}
\tilde h(t)=\bJ\varphi=\bJ\theta'=(1,-1,0,\ldots,0).
\end{equation}
Then, by Theorem~\ref{T:ie}, the function $\theta'(t)$ is the unique solution of the integral equation
\begin{equation}\label{e:rn-ie}
(\bI+\bN^\ast+\bJ)\theta'=\tilde h(t).
\end{equation}

In view of~(\ref{e:A}),~(\ref{e:thet-rho}) and~(\ref{e:rn-w}), it follows from~(\ref{e:rn-wb1}) that the boundary values of the function $F(z)$ satisfy 
\begin{equation}\label{e:rn-AF=}
A(t)F(\eta(t)) =\gamma(t)+h(t)+\i[(\rho(t)+\mu(t))+\nu(t)]
\end{equation} 
where 
\[
h(t)=\left(\ln\frac{1}{c},\ln\frac{R_1}{c},\ldots,\ln\frac{R_m}{c}\right)\in S,
\]
\[
\nu(t)=(c_0,c_1,\ldots,c_m)\in S,
\]
and the real functions $\gamma$, $\mu$ are defined by
\begin{equation}\label{e:rn-gam}
\gamma(t)+\i\mu(t):= -\log\left(1-\frac{\eta(t)}{z_1}\right).
\end{equation}

The boundary values of the mapping function $\omega$ are given by~(\ref{e:rn-w}). By~(\ref{e:rn-w'}), the boundary values of the derivative $\omega'$ are given by
\begin{equation}\label{e:rn-w'2}
\omega'(\eta(t))=\frac{\i\theta'(t)R(t)e^{\i\theta(t)}}{\dot\eta(t)}.
\end{equation}
Then, we have for $z\in G$,
\begin{equation}\label{e:rn-b-w(z)-1}
\omega(z)=\frac{1}{2\pi\i}\int_{\Gamma}\frac{\omega(\eta)}{\eta-z}d\eta =\frac{1}{2\pi\i}\int_{J}\frac{R(t)e^{\i\theta(t)}}{\eta(t)-z} \dot\eta(t)dt,
\end{equation}
and
\begin{equation}\label{e:rn-b-w'(z)-1}
\omega'(z)=\frac{1}{2\pi\i}\int_{\Gamma}\frac{\omega'(\eta)}{\eta-z}d\eta =\frac{1}{2\pi}\int_{J}\frac{\theta'(t)R(t)e^{\i\theta(t)}}{\eta(t)-z} dt.
\end{equation}

If $w\in\Omega$, then by the Cauchy integral formula,
\begin{equation}\label{e:rn-b-w(z)-i1}
\omega^{-1}(w)=\frac{1}{2\pi\i}\int_{\partial\Omega}\frac{\omega^{-1}(\xi)}{\xi-w}d\xi, 
\end{equation}
which on introducing $\xi=R(\theta)e^{\i\theta}$ and $\theta=\theta(t)$, we obtain
\begin{equation}\label{e:rn-b-w(z)-i3}
\omega^{-1}(w)= \frac{1}{2\pi\i}\int_{J} \frac{\eta(t)}{R(t)e^{\i\theta(t)}-w}R(t)e^{\i\theta(t)}\i\theta'(t)dt,
\end{equation} 
where we use the notation $R(\theta(t))=R(t)$ since $R$ is constant on each interval $J_j$, $j=0,1,\ldots,m$.

\section{A disc with circular slit region}

This canonical region $\Omega$ is the interior of the unit circle which has been
slit along $m$ arcs of circles. We assume that $\omega$ maps the curve $\Gamma_0$ onto the unit circle $|w|=1$ and the curves $\Gamma_j$, $j=1,2,\ldots,m$, onto slits on the circles $|w|=R_j$, where $R_1,\ldots,R_m$ are undetermined real constants. Then, the boundary values of the mapping function $\omega$ are given by
\begin{equation}\label{e:dc-w}
\omega(\eta(t))=R(t)e^{\i\theta(t)}
\end{equation}
where $\theta(t)$ is the boundary correspondence function and $R(t)=(1,R_1,\ldots,R_m)$. Thus, by differentiating both sides of~(\ref{e:dc-w}), we obtain
\begin{equation}\label{e:dc-w'}
\dot\eta(t)\frac{\omega'(\eta(t))}{\omega(\eta(t))}=\i\theta'(t).
\end{equation}

The function $\omega(z)$ is uniquely determined by assuming
\begin{equation}\label{e:dc-b-cd1}
\omega(0)=0, \quad \omega'(0)>0.
\end{equation}
Thus the function $\omega(z)$ can be expressed in the form
\begin{equation}\label{e:dc-wb}
\omega(z)= c\,z\,e^{zF(z)}
\end{equation}
where $c:=\omega'(0)$ is an undetermined positive real constant and $F(z)$ is an auxiliary unknown function. Hence
\begin{equation}\label{e:dc-wb1}
\eta(t)F(\eta(t))+\ln c+\log\eta(t) =\log(\omega(\eta(t))).
\end{equation} 
Thus, by differentiating both sides of~(\ref{e:dc-wb1}) and using~(\ref{e:dc-w'}), we obtain
\begin{equation}\label{e:dc-b-bd}
\dot\eta(t)\left(F(\eta(t))+\eta(t)F'(\eta(t))\right) +\frac{\dot\eta(t)}{\eta(t)} =\i\theta'(t).
\end{equation}

The function $f(z)$ defined in $G$ by
\begin{equation}\label{e:dc-f}
f(z):= zF(z)+z^2F'(z)+1, 
\end{equation}
and the function $g(z)$ defined in $G^-$ by $g(z):=0$ satisfy the assumptions of Theorem~\ref{T:ie} and their boundary values satisfy~(\ref{e:tAf+tAg}) with 
\begin{equation}\label{e:dc-b-gam-mu}
\upsilon(t)=0\quad{\rm and}\quad\varphi(t)=\theta'(t).
\end{equation}
Since the image of the curve $\Gamma_0$ is counterclockwise oriented and the images of the curves $\Gamma_j$, $j=1,\ldots,m$, are slits so we have $\theta_0(2\pi)-\theta_0(0)=2\pi$ and $\theta_j(2\pi)-\theta_j(0)=0$. Hence the function $\tilde h(t)$ in~(\ref{e:ie}) is given by
\begin{equation}\label{e:dc-b-th}
\tilde h(t)=\bJ\varphi=\bJ\theta'=(1,0,\ldots,0).
\end{equation}
Then, by Theorem~\ref{T:ie}, the function $\theta'(t)$ is the unique solution of the integral equation
\begin{equation}\label{e:dc-ie}
(\bI+\bN^\ast+\bJ)\theta'=\tilde h(t).
\end{equation}

In view of~(\ref{e:A}),~(\ref{e:thet-rho}) and~(\ref{e:dc-w}), it follows from~(\ref{e:dc-wb1}) that the boundary values of the function $F(z)$ satisfy 
\begin{equation}\label{e:dc-AF=}
A(t)F(\eta(t)) =\gamma(t)+h(t)+\i[(\rho(t)+\mu(t))+\nu(t)]
\end{equation} 
where 
\[
h(t)=\left(\ln\frac{1}{c},\ln\frac{R_1}{c},\ldots,\ln\frac{R_m}{c}\right)\in S,
\]
\[
\nu(t)=(c_0,c_1,\ldots,c_m)\in S,
\]
and the real functions $\gamma$, $\mu$ are defined by
\begin{equation}\label{e:dc-gam}
\gamma(t)+\i\mu(t):= -\log(\eta(t)).
\end{equation}

The boundary values of the mapping function $\omega$ are given by~(\ref{e:dc-w}). By~(\ref{e:dc-w'}), the boundary values of the derivative $\omega'$ are given by
\begin{equation}\label{e:dc-w'2}
\omega'(\eta(t))=\frac{\i\theta'(t)R(t)e^{\i\theta(t)}}{\dot\eta(t)}.
\end{equation}
Hence, the values of mapping function $\omega(z)$ and its derivative $\omega'(z)$ are given for $z\in G$ by
\begin{equation}\label{e:dc-b-w(z)-1}
\omega(z)=\frac{1}{2\pi\i}\int_{\Gamma}\frac{\omega(\eta)}{\eta-z}d\eta =\frac{1}{2\pi\i}\int_{J}\frac{R(t)e^{\i\theta(t)}}{\eta(t)-z} \dot\eta(t)dt,
\end{equation}
and
\begin{equation}\label{e:dc-b-w'(z)-1}
\omega'(z)=\frac{1}{2\pi\i}\int_{\Gamma}\frac{\omega'(\eta)}{\eta-z}d\eta =\frac{1}{2\pi}\int_{J}\frac{\theta'(t)R(t)e^{\i\theta(t)}}{\eta(t)-z} dt.
\end{equation}

The values of inverse mapping function $\omega^{-1}(w)$ are given for $w\in\Omega$ by
\begin{equation}\label{e:dc-b-w(z)-i1}
\omega^{-1}(w)=\frac{1}{2\pi\i}\int_{\partial\Omega} \frac{\omega^{-1}(\xi)}{\xi-w}d\xi, 
\end{equation}
which on introducing $\xi=R(\theta)e^{\i\theta}$ and $\theta=\theta(t)$ becomes
\begin{equation}\label{e:dc-b-w(z)-i3}
\omega^{-1}(w)= \frac{1}{2\pi\i}\int_{J} \frac{\eta(t)}{R(t)e^{\i\theta(t)}-w}R(t)e^{\i\theta(t)}\i\theta'(t)dt.
\end{equation}

\section{Circular slit region}

This canonical region $\Omega$ is the entire $w-$plane with $m+1$ slits along the circles $|w|=R_k$ where $R_0,R_1,\ldots,R_m$ are undetermined real constants. 
Then, the boundary values of the mapping function $\omega$ are given by
\begin{equation}\label{e:cr-w}
\omega(\eta(t))=R(t)e^{\i\theta(t)}
\end{equation}
where $\theta(t)$ is the boundary correspondence function and $R(t)=(R_0,R_1,\ldots,R_m)$. Thus, by differentiating both sides of~(\ref{e:cr-w}), we obtain
\begin{equation}\label{e:cr-w'}
\dot\eta(t)\frac{\omega'(\eta(t))}{\omega(\eta(t))}=\i\theta'(t).
\end{equation}

The function $\omega$ can be uniquely determined by assuming
\begin{equation}\label{e:cr-b-cd}
\omega(\alpha) = 0, \quad \omega(0)=\infty, \quad
\Res_{z=0} \omega(z) = 1,
\end{equation}
where $\alpha$ is a fixed point in $G$. Hence $\omega$ can be written as
\begin{equation}\label{e:cr-wb}
\omega(z) = \left(\frac{1}{z}-\frac{1}{\alpha}\right)\,
e^{z F(z)}
\end{equation}
where $F(z)$ is an auxiliary unknown function. Hence
\begin{equation}\label{e:cr-wb1}
\eta(t)F(\eta(t))+\log\left(\frac{1}{\eta(t)}-\frac{1}{\alpha}\right) =\log(\omega(\eta(t))).
\end{equation} 
Thus, by differentiating both sides of~(\ref{e:cr-wb1}) and using~(\ref{e:cr-w'}), we obtain
\begin{equation}\label{e:cr-b-bd}
\dot\eta(t)\left(F(\eta(t))+\eta(t)F'(\eta(t))\right) +\frac{\dot\eta(t)}{\eta(t)}\frac{\alpha}{\eta(t)-\alpha} =\i\theta'(t).
\end{equation}

The function $f(z)$ defined in $G$ by
\begin{equation}\label{e:cr-f}
f(z):= zF(z)+z^2F'(z), 
\end{equation}
and the function $g(z)$ defined in $G^-$ by
\begin{equation}\label{e:cr-g}
g(z):= \frac{\alpha}{z-\alpha},
\end{equation}
satisfy the assumptions of Theorem~\ref{T:ie} and their boundary values satisfy~(\ref{e:tAf+tAg}) with 
\begin{equation}\label{e:cr-b-gam-mu}
\upsilon(t)=0\quad{\rm and}\quad\varphi(t)=\theta'(t).
\end{equation}
Since the image of the curves $\Gamma_j$, $j=0,1,\ldots,m$, are slits so we have $\theta_j(2\pi)-\theta_j(0)=0$. Hence the function $\tilde h(t)$ in~(\ref{e:ie}) is given by
\begin{equation}\label{e:cr-b-th}
\tilde h(t)=\bJ\varphi=\bJ\theta'=(0,0,\ldots,0).
\end{equation}
Then, by Theorem~\ref{T:ie}, the function $\theta'(t)$ is the unique solution of the integral equation
\begin{equation}\label{e:cr-ie}
(\bI+\bN^\ast+\bJ)\theta'=2\phi
\end{equation}
where
\[
\phi(t) :=\Im[\tilde A(t)g(\eta(t))] =\Im\left[\frac{\dot\eta(t)}{\eta(t)}\frac{\alpha}{\eta(t)-\alpha}\right].
\]

In view of~(\ref{e:A}),~(\ref{e:thet-rho}) and~(\ref{e:cr-w}), it follows from~(\ref{e:cr-wb1}) that the boundary values of the function $F(z)$ satisfy 
\begin{equation}\label{e:cr-AF=}
A(t)F(\eta(t)) =\gamma(t)+h(t)+\i[(\rho(t)+\mu(t))+\nu(t)]
\end{equation} 
where 
\[
h(t)=(\ln R_0,\ln R_1,\ldots,\ln R_m)\in S,
\]
\[
\nu(t)=(c_0,c_1,\ldots,c_m)\in S,
\]
and the real functions $\gamma$, $\mu$ are defined by
\begin{equation}\label{e:cr-gam}
\gamma(t)+\i\mu(t):= -\log\left(\frac{1}{\eta(t)}-\frac{1}{\alpha}\right). 
\end{equation}

The boundary values of the mapping function $\omega$ are given by~(\ref{e:cr-w}). By~(\ref{e:cr-w'}), the boundary values of the derivative $\omega'$ are given by
\begin{equation}\label{e:cr-w'2}
\omega'(\eta(t))=\frac{\i\theta'(t)R(t)e^{\i\theta(t)}}{\dot\eta(t)}.
\end{equation}
Since $\Res_{z=0} \omega(z)=1$, then, by~(\ref{e:cau-b1}) and~(\ref{e:cau-b1'}), the values of the mapping function $\omega(z)$ and its derivative $\omega'(z)$ are given for $z\in G$ by
\begin{equation}\label{e:cr-b-w(z)-1}
\omega(z)=\frac{1}{z}+\frac{1}{2\pi\i}\int_{\Gamma}\frac{\omega(\eta)}{\eta-z}d\eta =\frac{1}{z}+\frac{1}{2\pi\i}\int_{J}\frac{R(t)e^{\i\theta(t)}}{\eta(t)-z} \dot\eta(t)dt,
\end{equation}
and
\begin{equation}\label{e:cr-b-w'(z)-1}
\omega'(z)=-\frac{1}{z^2}+\frac{1}{2\pi\i}\int_{\Gamma}\frac{\omega'(\eta)}{\eta-z}d\eta =-\frac{1}{z^2}+\frac{1}{2\pi}\int_{J}\frac{\theta'(t)R(t)e^{\i\theta(t)}}{\eta(t)-z} dt.
\end{equation}

For $w\in\Omega$, then by the Cauchy integral formula and since $\omega^{-1}(\infty)=0$, the values of the inverse mapping function $\omega^{-1}(w)$ are given by
\begin{equation}\label{e:cr-b-w(z)-i1}
\omega^{-1}(w)=\frac{1}{2\pi\i}\int_{\partial\Omega}\frac{\omega^{-1}(\xi)}{\xi-w}d\xi, 
\end{equation}
which on introducing $\xi=R(\theta)e^{\i\theta}$ and $\theta=\theta(t)$ becomes
\begin{equation}\label{e:cr-b-w(z)-i3}
\omega^{-1}(w)= \frac{1}{2\pi\i}\int_{J} \frac{\eta(t)}{R(t)e^{\i\theta(t)}-w} R(t)e^{\i\theta(t)}\i\theta'(t)dt.
\end{equation}

\section{Radial slit region}

This canonical region $\Omega$ is the entire $w-$plane with $m+1$ slits along the rays $\arg(w)=R_k$ where $R_k$, $k=0,1,\ldots, m$, are undetermined real constants. 
Then, the boundary values of the mapping function $\omega$ are given by
\begin{equation}\label{e:rd-w}
\omega(\eta(t))=e^{\theta(t)}e^{\i R(t)}
\end{equation}
where $\theta(t)$ is the boundary correspondence function and $R(t)=(R_0,R_1,\ldots,R_m)$. We use the notation $e^{\theta(t)}$ to emphasize that the functions multiplied by $e^{\i R_k}$ are positive. Thus, by differentiating both sides of~(\ref{e:rd-w}), we obtain
\begin{equation}\label{e:rd-w'}
\dot\eta(t)\frac{\omega'(\eta(t))}{\omega(\eta(t))} =\theta'(t).
\end{equation}

The function $\omega$ can be uniquely determined by assuming
\begin{equation}\label{e:rd-b-cd}
\omega(\alpha) = 0, \quad \omega(0)=\infty, \quad
\Res_{z=0} \omega(z) = 1,
\end{equation}
where $\alpha$ is a fixed point in $G$. Hence $\omega$ can be written as
\begin{equation}\label{e:rd-wb}
\omega(z) = \left(\frac{1}{z}-\frac{1}{\alpha}\right)\,
e^{\i zF(z)}
\end{equation}
where $F(z)$ is an auxiliary unknown function. Hence
\begin{equation}\label{e:rd-wb1}
-\eta(t)F(\eta(t))+\i\log\left(\frac{1}{\eta(t)}-\frac{1}{\alpha}\right) =\i\log(\omega(\eta(t))).
\end{equation} 
Thus, by differentiating both sides of~(\ref{e:rd-wb1}) and using~(\ref{e:rd-w'}), we obtain
\begin{equation}\label{e:rd-b-bd}
\dot\eta(t)\left(-F(\eta(t))-\eta(t)F'(\eta(t))\right) +\frac{\dot\eta(t)}{\eta(t)}\frac{\alpha\i}{\eta(t)-\alpha}  =\i\theta'(t).
\end{equation}

The function $f(z)$ defined in $G$ by
\begin{equation}\label{e:rd-f}
f(z):= -zF(z)-z^2F'(z),
\end{equation}
and the function $g(z)$ defined in $G^-$ by
\begin{equation}\label{e:rd-g}
g(z):=  \frac{\alpha\i}{z-\alpha},
\end{equation}
satisfy the assumptions of Theorem~\ref{T:ie} and their boundary values satisfy~(\ref{e:tAf+tAg}) with 
\begin{equation}\label{e:rd-b-gam-mu}
\upsilon(t)=0\quad{\rm and}\quad\varphi(t)=\theta'(t).
\end{equation}
Since the image of the curves $\Gamma_j$, $j=0,1,\ldots,m$, are slits so we have $\theta_j(2\pi)-\theta_j(0)=0$. Hence the function $\tilde h(t)$ in~(\ref{e:ie}) is given by
\begin{equation}\label{e:rd-b-th}
\tilde h(t)=\bJ\varphi=\bJ\theta' =(0,0,\ldots,0).
\end{equation}
Then, by Theorem~\ref{T:ie}, the function $\theta'(t)$ is the unique solution of the integral equation
\begin{equation}\label{e:rd-ie}
(\bI+\bN^\ast+\bJ)\theta' =2\phi
\end{equation}
where 
\[
\phi(t) :=\Im[\tilde A(t)g(\eta(t))] =\Im\left[\frac{\dot\eta(t)}{\eta(t)}\frac{\alpha\i}{\eta(t)-\alpha}\right].
\]

It follows from~(\ref{e:A}),~(\ref{e:rd-w}) and~(\ref{e:rd-wb1}) that the boundary values of the function $F(z)$ satisfy 
\begin{equation}
A(t)F(\eta(t)) =\i\log\left(\frac{1}{\eta(t)}-\frac{1}{\alpha}\right)-\i(\theta(t)+\i R(t)).
\end{equation} 
Thus, in view of~(\ref{e:thet-rho}), the boundary values of the function $F(z)$ are given by 
\begin{equation}\label{e:rd-AF=}
A(t)F(\eta(t)) =\gamma(t)+h(t)+\i[(-\rho(t)+\mu(t))+\nu(t)]
\end{equation} 
where 
\[
h(t)=(R_0,R_1,\ldots,R_m)\in S,
\]
\[
\nu(t)=-(c_0,c_1,\ldots,c_m)\in S,
\]
and the real functions $\gamma$, $\mu$ are defined by
\begin{equation}\label{e:rd-gam}
\gamma(t)+\i\mu(t):= \i\log\left(\frac{1}{\eta(t)}-\frac{1}{\alpha}\right).
\end{equation}

The boundary values of the mapping function $\omega$ are given by~(\ref{e:rd-w}). By~(\ref{e:rd-w'}), the boundary values of the derivative $\omega'$ are given by
\begin{equation}\label{e:rd-w'2}
\omega'(\eta(t))=\frac{\theta'(t)e^{\theta(t)}e^{\i R(t)}}{\dot\eta(t)}.
\end{equation}
Since $\Res_{z=0} \omega(z)=1$, then, by~(\ref{e:cau-b1}) and~(\ref{e:cau-b1'}), the values of the mapping function $\omega(z)$ and its derivative $\omega'(z)$ are given for $z\in G$ by
\begin{equation}\label{e:rd-b-w(z)-1}
\omega(z)=\frac{1}{z}+\frac{1}{2\pi\i}\int_{\Gamma}\frac{\omega(\eta)}{\eta-z}d\eta =\frac{1}{z}+\frac{1}{2\pi\i}\int_{J}\frac{e^{\theta(t)} e^{\i R(t)}}{\eta(t)-z} \dot\eta(t)dt,
\end{equation}
and
\begin{equation}\label{e:rd-b-w'(z)-1}
\omega'(z)=-\frac{1}{z^2}+\frac{1}{2\pi\i}\int_{\Gamma}\frac{\omega'(\eta)}{\eta-z}d\eta =-\frac{1}{z^2}+\frac{1}{2\pi\i} \int_{J}\frac{\theta'(t)e^{\theta(t)} e^{\i R(t)}}{\eta(t)-z} dt.
\end{equation}

For $w\in\Omega$, then by the Cauchy integral formula and since $\omega^{-1}(\infty)=0$, the values of the inverse mapping function $\omega^{-1}(w)$ are given by
\begin{equation}\label{e:rd-b-w(z)-i1}
\omega^{-1}(w)=\frac{1}{2\pi\i}\int_{\partial\Omega}\frac{\omega^{-1}(\xi)}{\xi-w}d\xi, 
\end{equation}
which on introducing $\xi=e^{\theta}e^{\i R(\theta)}$ and $\theta=\theta(t)$ becomes
\begin{equation}\label{e:rd-b-w(z)-i3}
\omega^{-1}(w)= \frac{1}{2\pi\i}\int_{J} \frac{\eta(t)}{e^{\theta(t)}e^{\i R(t)}-w} e^{\theta(t)}e^{\i R(t)}\theta'(t)dt.
\end{equation}

\section{Parallel slit region}
\label{sc:par}

This canonical region $\Omega$ is the entire $w-$plane with $m+1$ straight slits on the straight lines
\begin{equation}\label{e:par-lines}
\Re\left[e^{\i(\pi/2-\delta)}\, w\right]=R_j, \quad j=0,1,\ldots,m,
\end{equation}
where $R_0,R_1,\ldots,R_m$ are undetermined real constants and $\delta$ is the angle of intersection between the lines~(\ref{e:par-lines}) and the real axis. Thus, the boundary values of the mapping function $\omega$ satisfy
\begin{equation}\label{e:par-w}
e^{\i(\pi/2-\delta)} \omega(\eta(t)) = R(t)+\i\theta(t)
\end{equation}
where $\theta(t)$ is the boundary correspondence function and $R(t)=(R_0,R_1,\ldots,R_m)$. Thus, by differentiating both sides of~(\ref{e:par-w}), we obtain
\begin{equation}\label{e:par-w'}
e^{\i(\pi/2-\delta)}\dot\eta(t)\omega'(\eta(t))= \i\theta'(t).
\end{equation}

The mapping function $\omega$ is uniquely determined by the normalization
\begin{equation}\label{e:par-cd-b}
\omega(0) = \infty, \quad \lim_{z\to 0}\left(\omega(z)-1/z\right) =0.
\end{equation}
Thus, the function $\omega$ can be written as
\begin{equation}\label{e:par-wb}
\omega(z)=\frac{1}{z}+e^{-\i(\pi/2-\delta)}zF(z).
\end{equation}
where $F(z)$ is an auxiliary unknown function. Hence
\begin{equation}\label{e:par-wb1}
\eta(t)F(\eta(t))+\frac{e^{\i(\pi/2-\delta)}}{\eta(t)} =e^{\i(\pi/2-\delta)}\omega(\eta(t)).
\end{equation} 
Thus, by differentiating both sides of~(\ref{e:par-wb1}) and using~(\ref{e:par-w'}), we obtain
\begin{equation}\label{e:par-b-bd}
\dot\eta(t)\left(F(\eta(t))+\eta(t)F'(\eta(t))\right) +\dot\eta(t)\frac{-e^{\i(\pi/2-\delta)}}{\eta(t)^2}  =\i\theta'(t).
\end{equation}

The function $f(z)$ defined in $G$ by
\begin{equation}\label{e:par-f}
f(z):= zF(z)+z^2F'(z),
\end{equation}
and the function $g(z)$ defined in $G^-$ by
\begin{equation}\label{e:par-g}
g(z):= -\frac{e^{\i(\pi/2-\delta)}}{z}, 
\end{equation}
satisfy the assumptions of Theorem~\ref{T:ie} and their boundary values satisfy~(\ref{e:tAf+tAg}) with 
\begin{equation}\label{e:par-b-gam-mu}
\upsilon(t)=0\quad{\rm and}\quad\varphi(t)=\theta'(t).
\end{equation}
Since the image of the curves $\Gamma_j$, $j=0,1,\ldots,m$, are slits so we have $\theta_j(2\pi)-\theta_j(0)=0$. Hence the function $\tilde h(t)$ in~(\ref{e:ie}) is given by
\begin{equation}\label{e:par-b-th}
\tilde h(t)=\bJ\varphi=\bJ\theta'=(0,0,\ldots,0).
\end{equation}
Then, by Theorem~\ref{T:ie}, the function $\theta'$ is the unique solution of the integral equation
\begin{equation}\label{e:par-ie}
(\bI+\bN^\ast+\bJ)\theta'=2\phi
\end{equation}
where
\[
\phi(t) :=\Im[\tilde A(t)g(\eta(t))] =-\Im\left[\frac{\dot\eta(t)}{\eta(t)} \frac{e^{\i(\pi/2-\delta)}}{\eta(t)}\right].
\]

In view of~(\ref{e:A}),~(\ref{e:thet-rho}) and~(\ref{e:par-w}), it follows from~(\ref{e:par-wb1}) that the boundary values of the function $F(z)$ satisfy 
\begin{equation}\label{e:par-AF=}
A(t)F(\eta(t)) =\gamma(t)+h(t)+\i[(\rho(t)+\mu(t))+\nu(t)]
\end{equation} 
where 
\[
h(t)=(R_0,R_1,\ldots,R_m)\in S,
\]
\[
\nu(t)=(c_0,c_1,\ldots,c_m)\in S,
\]
and the real functions $\gamma$, $\mu$ are defined by
\begin{equation}\label{e:par-gam}
\gamma(t)+\i\mu(t):= -\frac{e^{\i(\pi/2-\delta)}}{\eta(t)}.
\end{equation}

The boundary values of the mapping function $\omega$ and its derivative $\omega'$ are given by~(\ref{e:rd-w}) and~(\ref{e:rd-w'}). Since $\Res_{z=0} \omega(z)=1$, then, by~(\ref{e:cau-b1}) and~(\ref{e:cau-b1'}), the values of the mapping function $\omega(z)$ and its derivative $\omega'(z)$ are given for $z\in G$ by
\begin{equation}\label{e:par-b-w(z)-1}
\omega(z)=\frac{1}{z}+\frac{1}{2\pi\i}\int_{\Gamma}\frac{\omega(\eta)}{\eta-z}d\eta =\frac{1}{z}+\frac{1}{2\pi\i}\int_{J}\frac{e^{-\i(\pi/2-\sigma)} [R(t)+\i\theta(t)]}{\eta(t)-z} \dot\eta(t)dt,
\end{equation}
and
\begin{equation}\label{e:par-b-w'(z)-1}
\omega'(z)=-\frac{1}{z^2}+\frac{1}{2\pi\i}\int_{\Gamma}\frac{\omega'(\eta)}{\eta-z}d\eta 
=-\frac{1}{z^2}+\frac{1}{2\pi}\int_{J}\frac{e^{-\i(\pi/2-\sigma)}\theta'(t)}{\eta(t)-z} dt.
\end{equation}

For $w\in\Omega$, then by the Cauchy integral formula and since $\omega^{-1}(\infty)=0$, the values of the inverse mapping function $\omega^{-1}(w)$ are given by
\begin{equation}\label{e:par-b-w(z)-i1}
\omega^{-1}(w)=\frac{1}{2\pi\i}\int_{\partial\Omega}\frac{\omega^{-1}(\xi)}{\xi-w}d\xi, 
\end{equation}
which on introducing $\xi=e^{-\i(\pi/2-\sigma)}[R(\theta)+\i\theta]$ and $\theta=\theta(t)$ becomes
\begin{equation}\label{e:par-b-w(z)-i3}
\omega^{-1}(w)= \frac{1}{2\pi\i}\int_{J} \frac{\eta(t)}{e^{-\i(\pi/2-\sigma)}[R(\theta)+\i\theta]-w} e^{-\i(\pi/2-\sigma)}\i\theta'(t)dt.
\end{equation}

\section{Numerical Example}
\label{sc:exp}

Since the functions $A_j$ and $\eta_j$ are $2\pi$-periodic, a reliable procedure for solving the integral equation~(\ref{e:ie}) numerically is by using the Nystr\"om method with the trapezoidal rule~\cite{atk}. The computational details are similar to previous works~\cite{nas-fun,nas-siam,nas-jmaa}. 

By using the trapezoidal rule with $n$ (an even positive integer) equidistant collocation points on each boundary component, solving the integral equation~(\ref{e:ie}) reduces to solving an $(m+1)n$ by $(m+1)n$ linear system. Since the integral equation~(\ref{e:ie}) is uniquely solvable, then for sufficiently large values of $n$ the obtained linear system is also uniquely solvable~\cite{atk}. In this paper, the linear system is solved using the Gauss elimination method. For calculating the function $\rho_k(t)$ in~(\ref{e:rho-F}), we approximate the function $\theta'_k(t)$ in~(\ref{e:thet'-F}) by the interpolating trigonometric polynomial of degree $n/2$ that interpolate $\theta'(t)$ at the $n$ equidistant points (see~\cite[p.~364]{wegc}). Then $\rho_k(t)$ is calculated by integrating the obtained interpolating trigonometric polynomial.

\begin{figure}%
\centerline{%
\hfill\scalebox{0.25}{\includegraphics{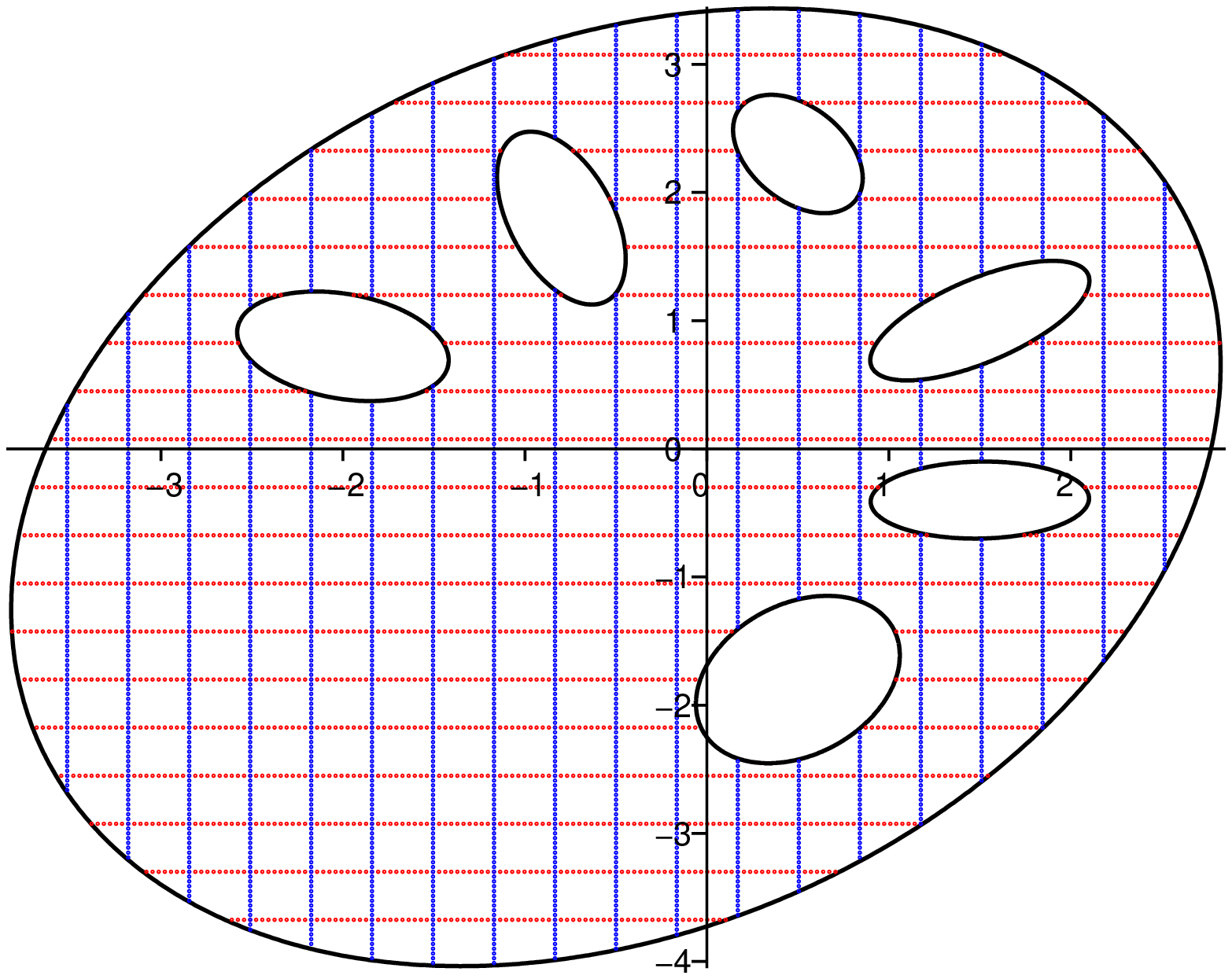}}%
\hfill\scalebox{0.25}{\includegraphics{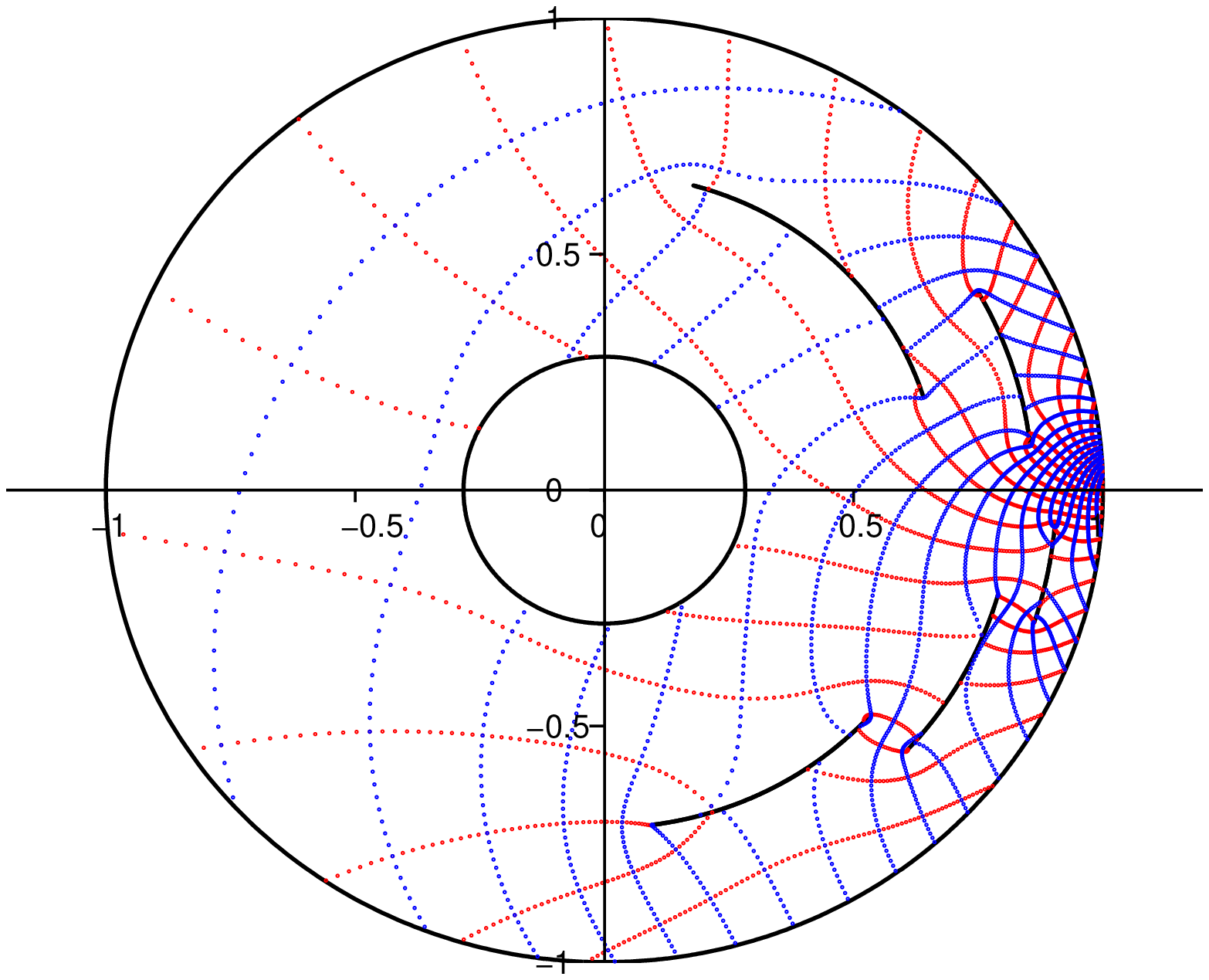}}%
\hfill\scalebox{0.25}{\includegraphics{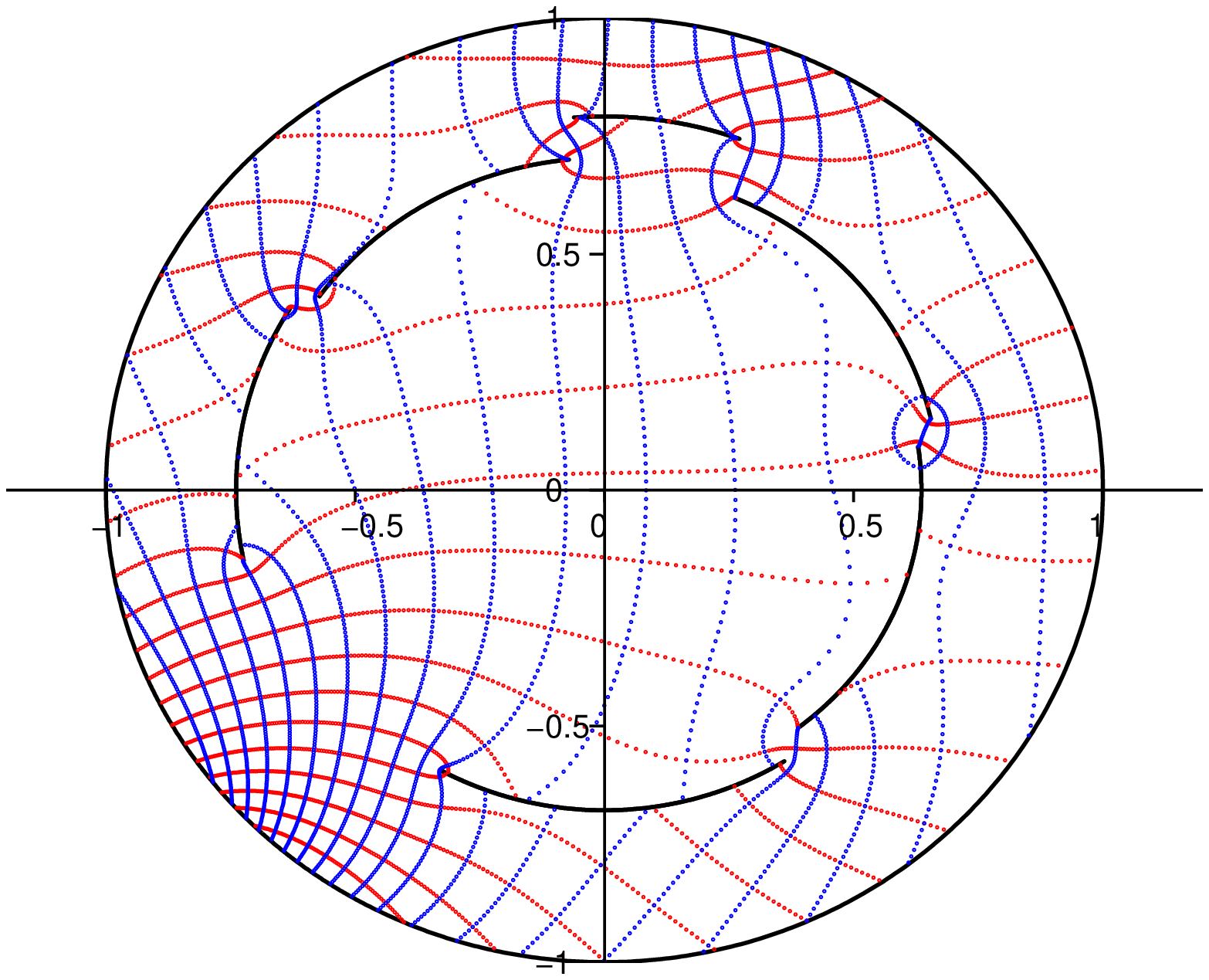}}%
\hfill}%
\centerline{%
\hfill\scalebox{0.25}{\includegraphics{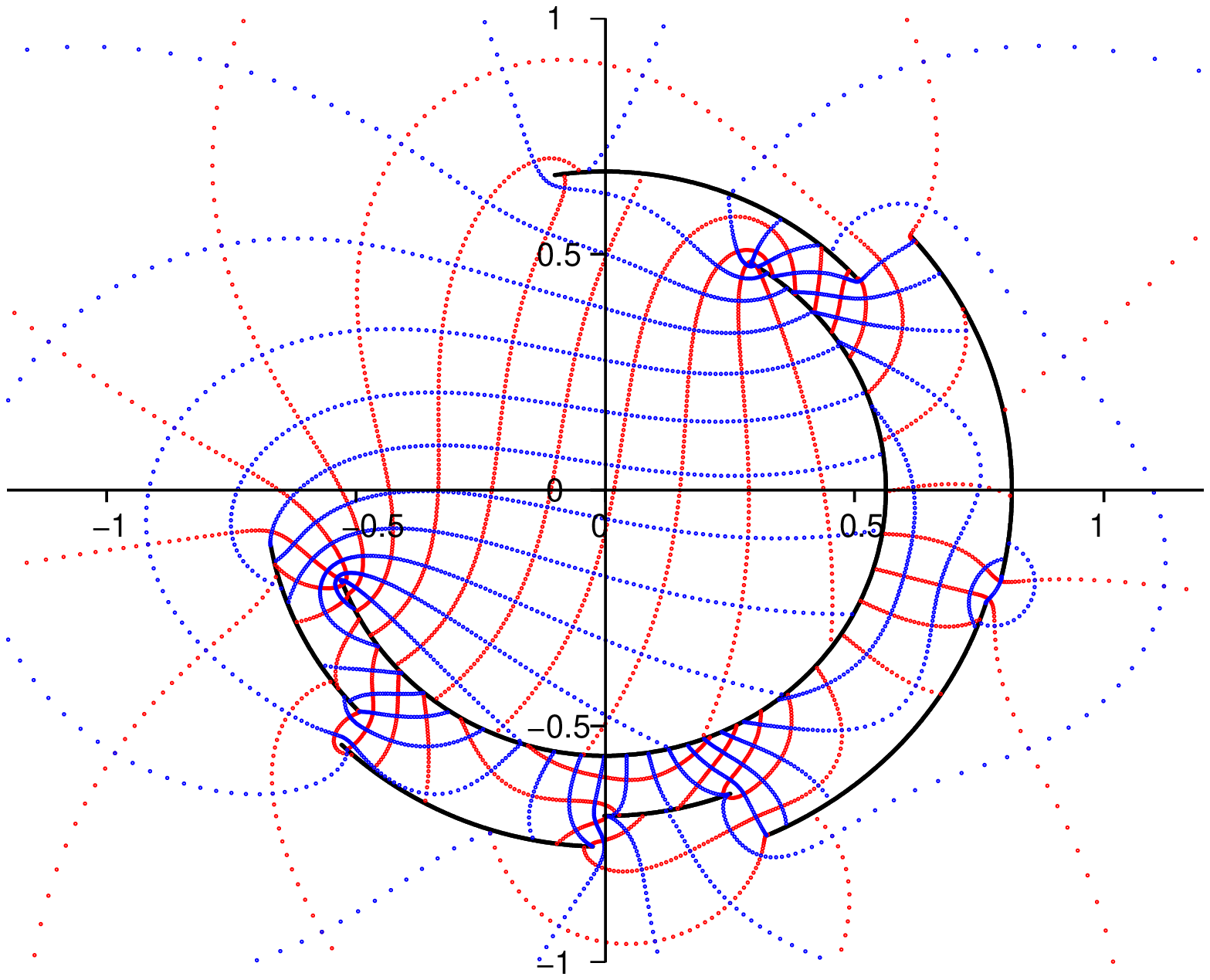}}%
\hfill\scalebox{0.25}{\includegraphics{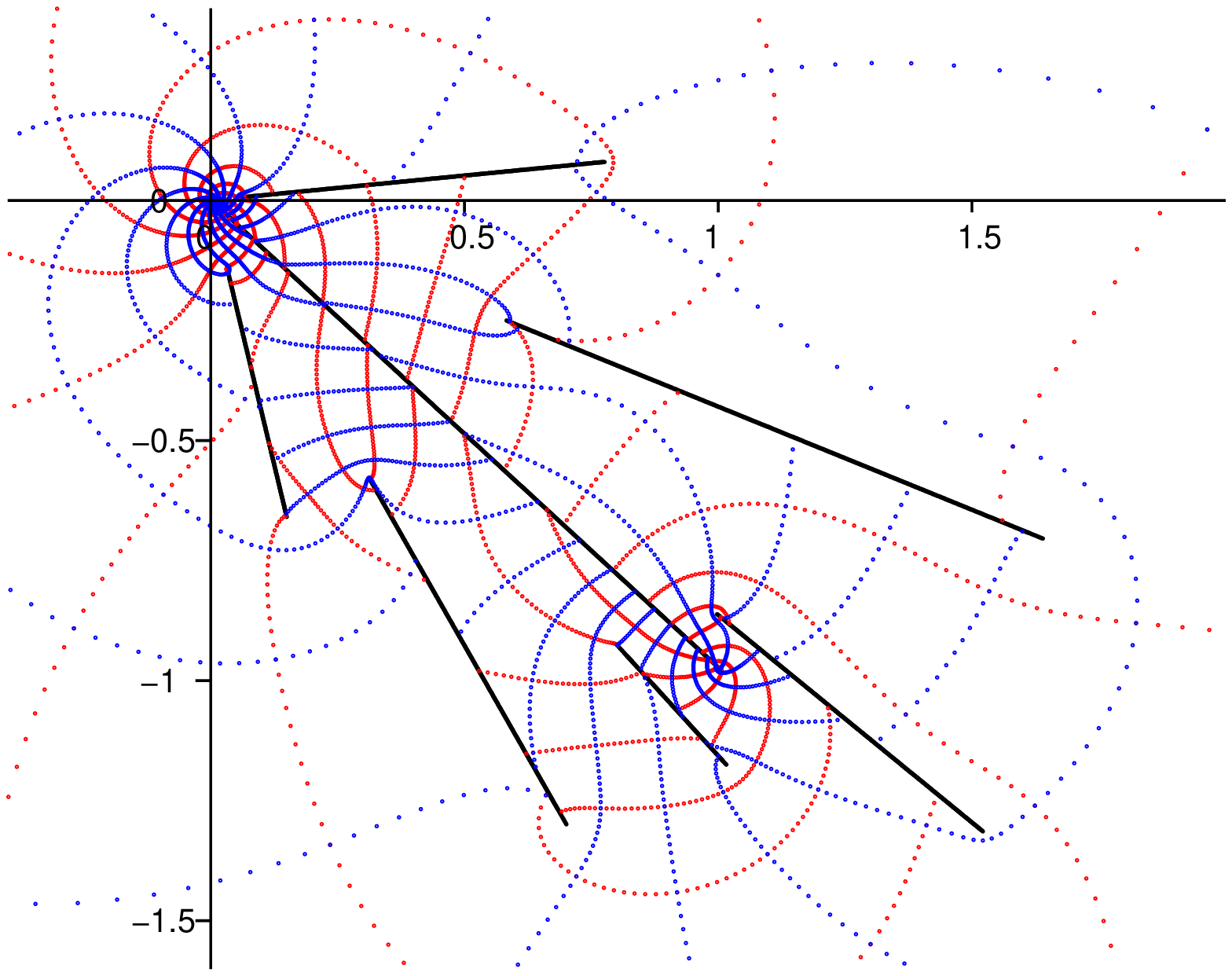}}%
\hfill\scalebox{0.25}{\includegraphics{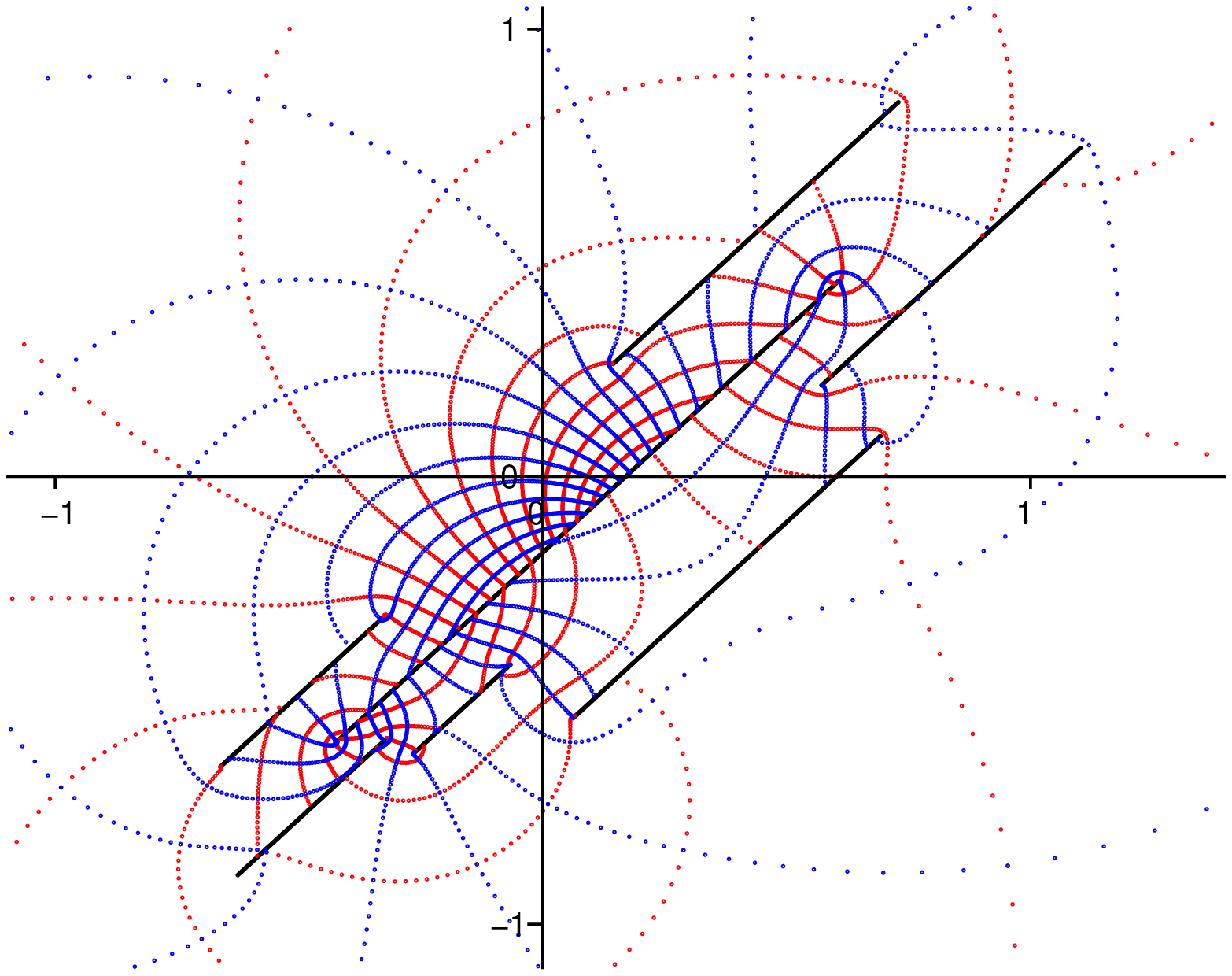}}%
\hfill}%
\caption{The original region $G$ and its images.} 
\label{f:e2or}
\end{figure}

\begin{figure}%
\centerline{%
\hfill\scalebox{0.25}{\includegraphics{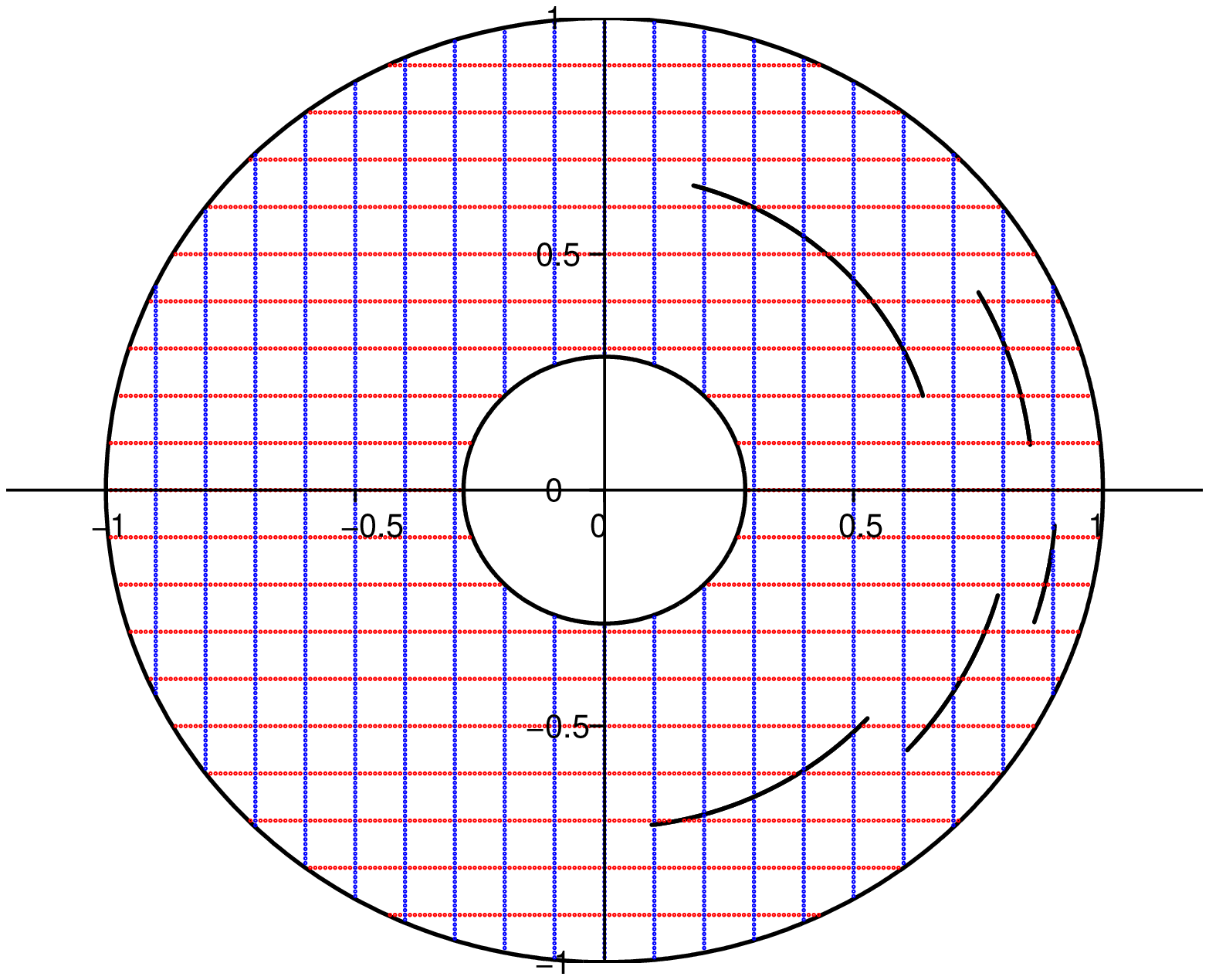}}%
\hfill\scalebox{0.25}{\includegraphics{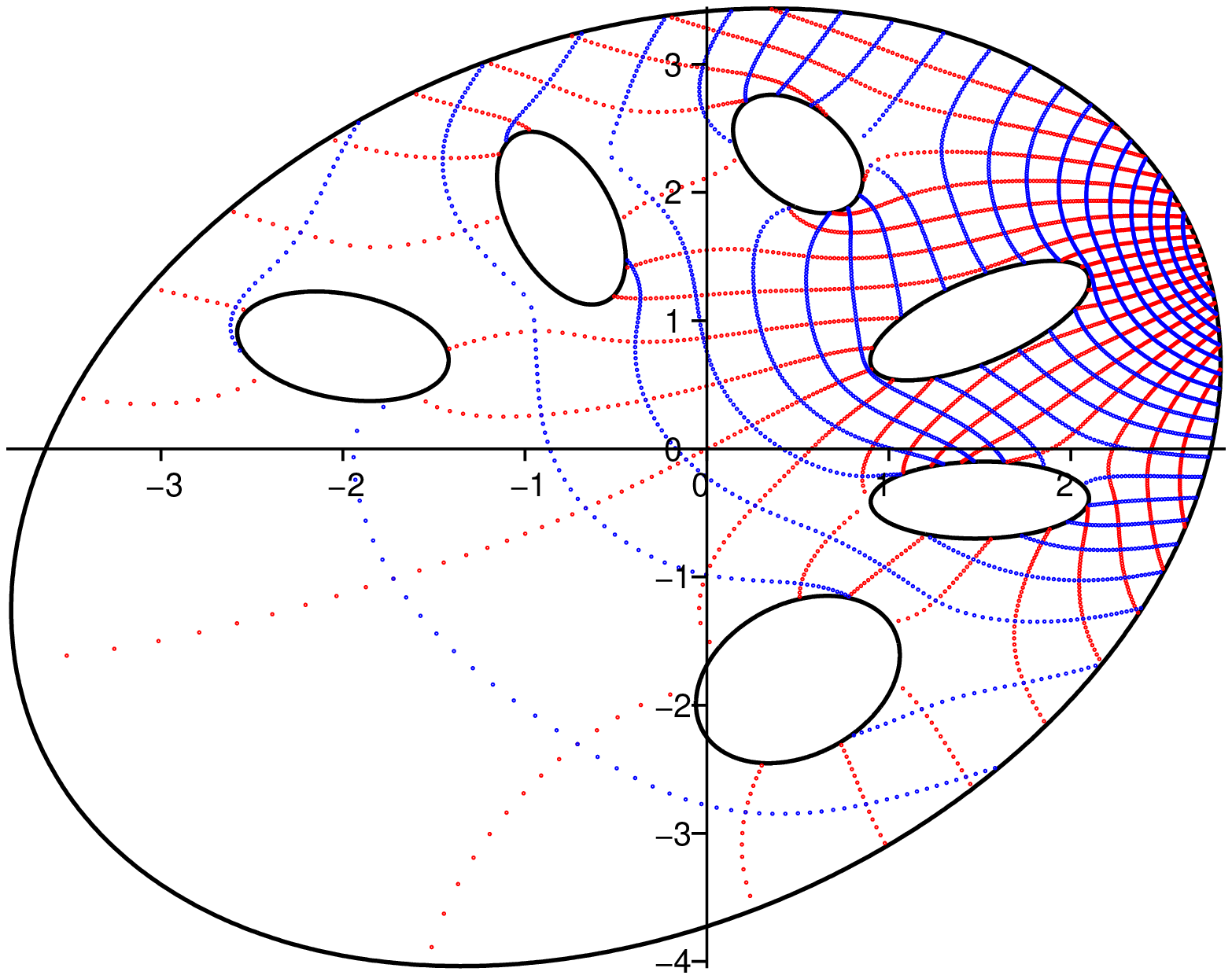}}%
\hfill}%
\caption{The inverse image of the annulus with circular slit canonical region.} 
\label{f:e2an}
\end{figure}

\begin{figure}%
\centerline{%
\hfill\scalebox{0.25}{\includegraphics{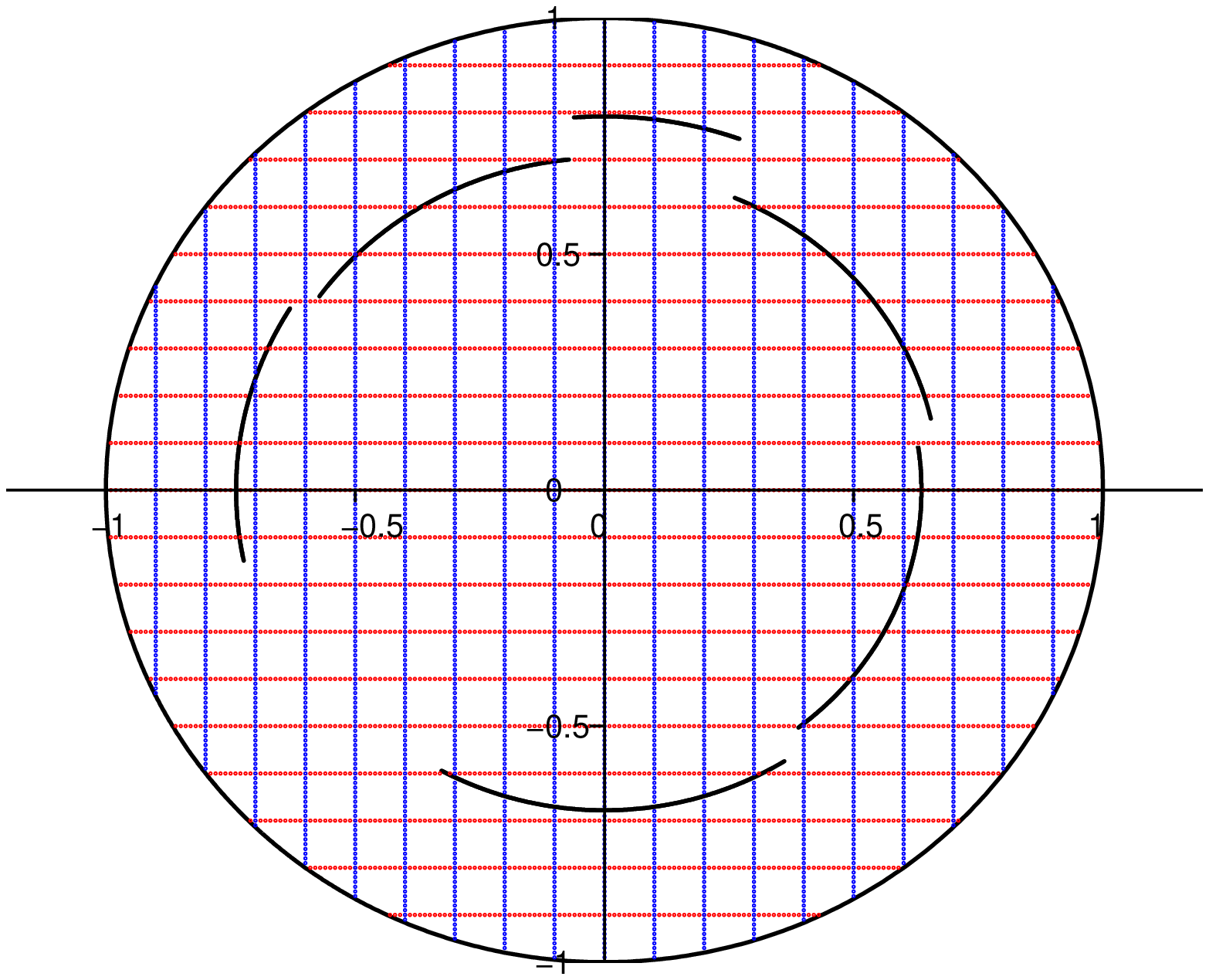}}%
\hfill\scalebox{0.25}{\includegraphics{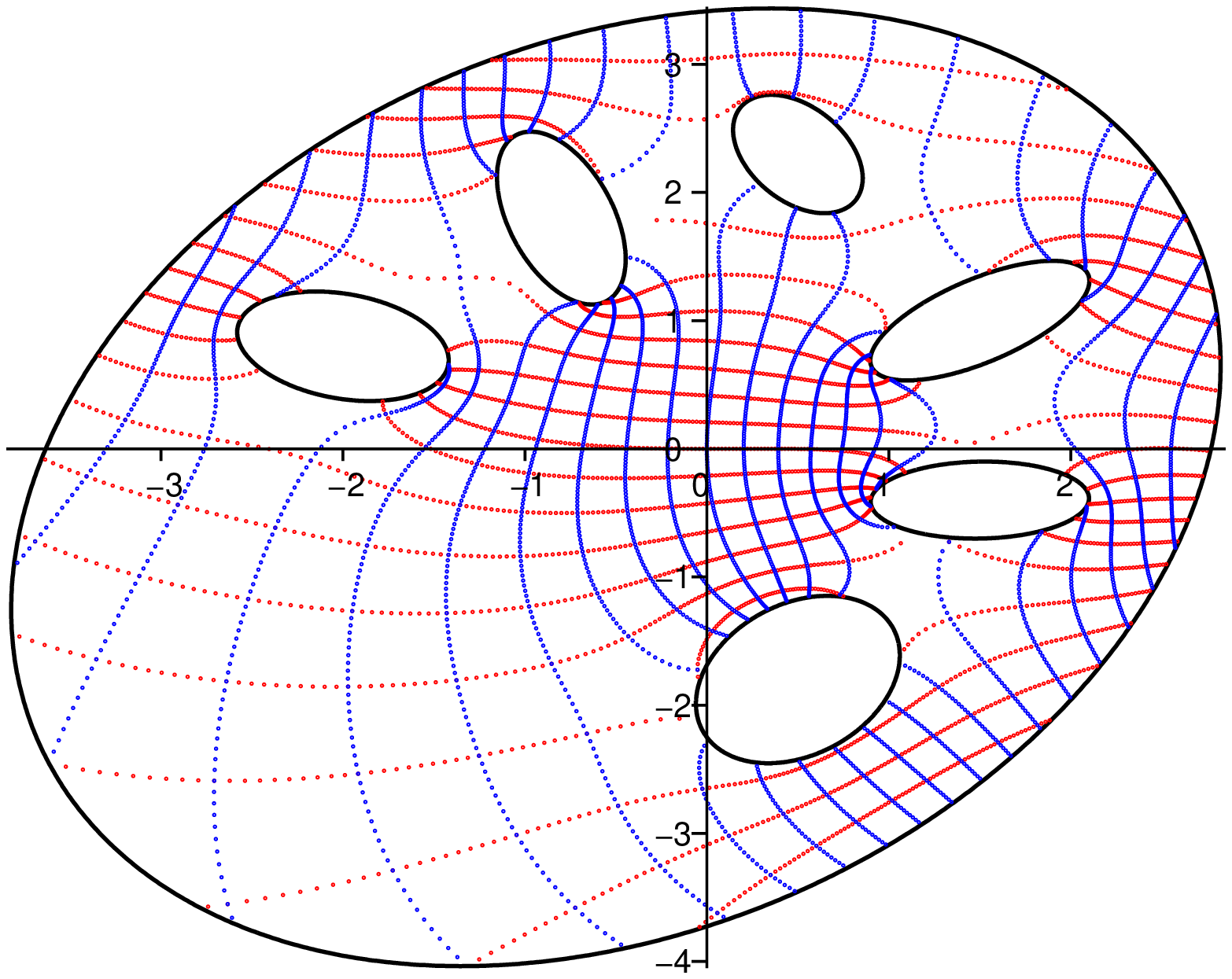}}%
\hfill}%
\caption{The inverse image of the disc with circular slit canonical region.} 
\label{f:e2dc}
\end{figure}

\begin{figure}%
\centerline{%
\hfill\scalebox{0.25}{\includegraphics{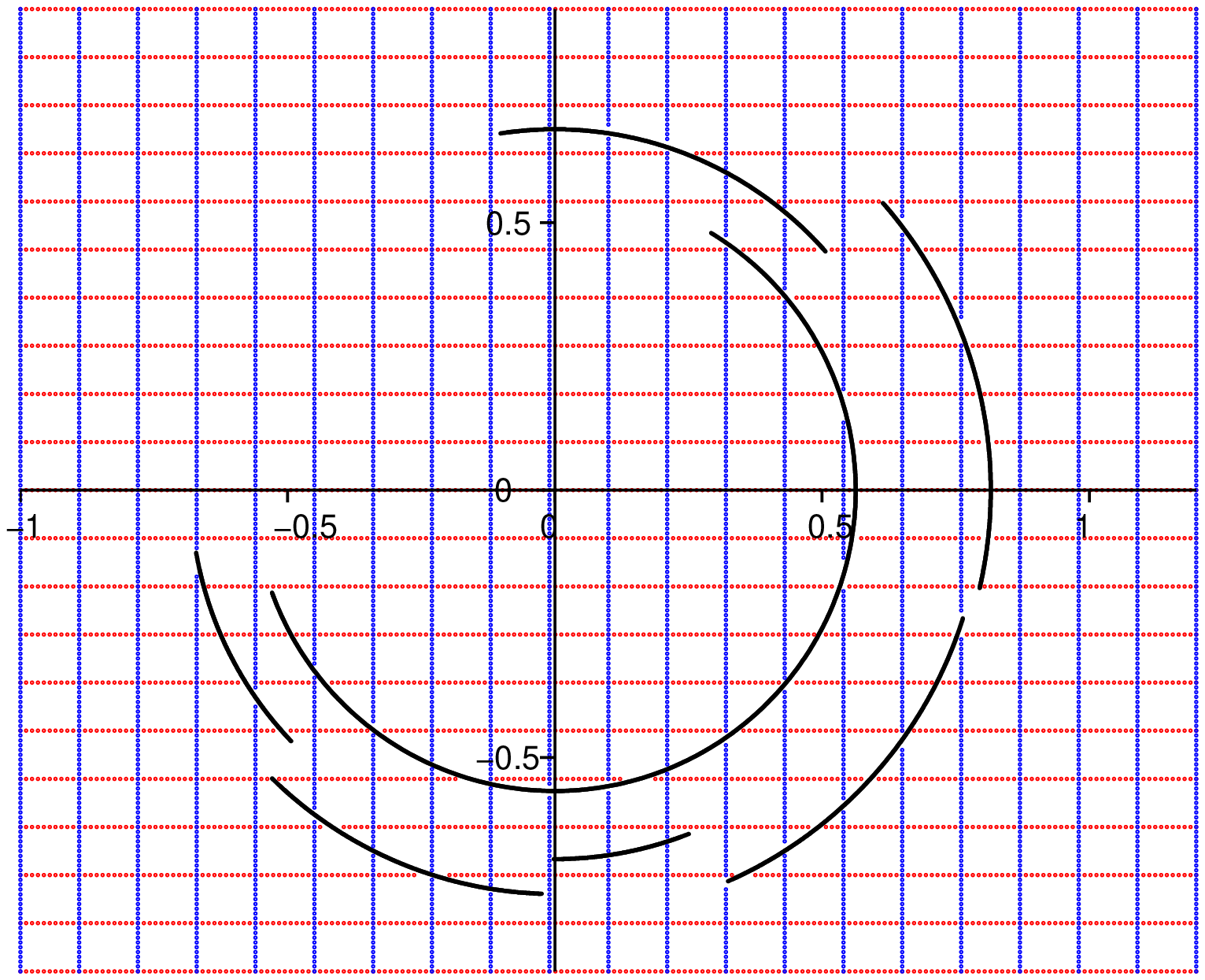}}%
\hfill\scalebox{0.25}{\includegraphics{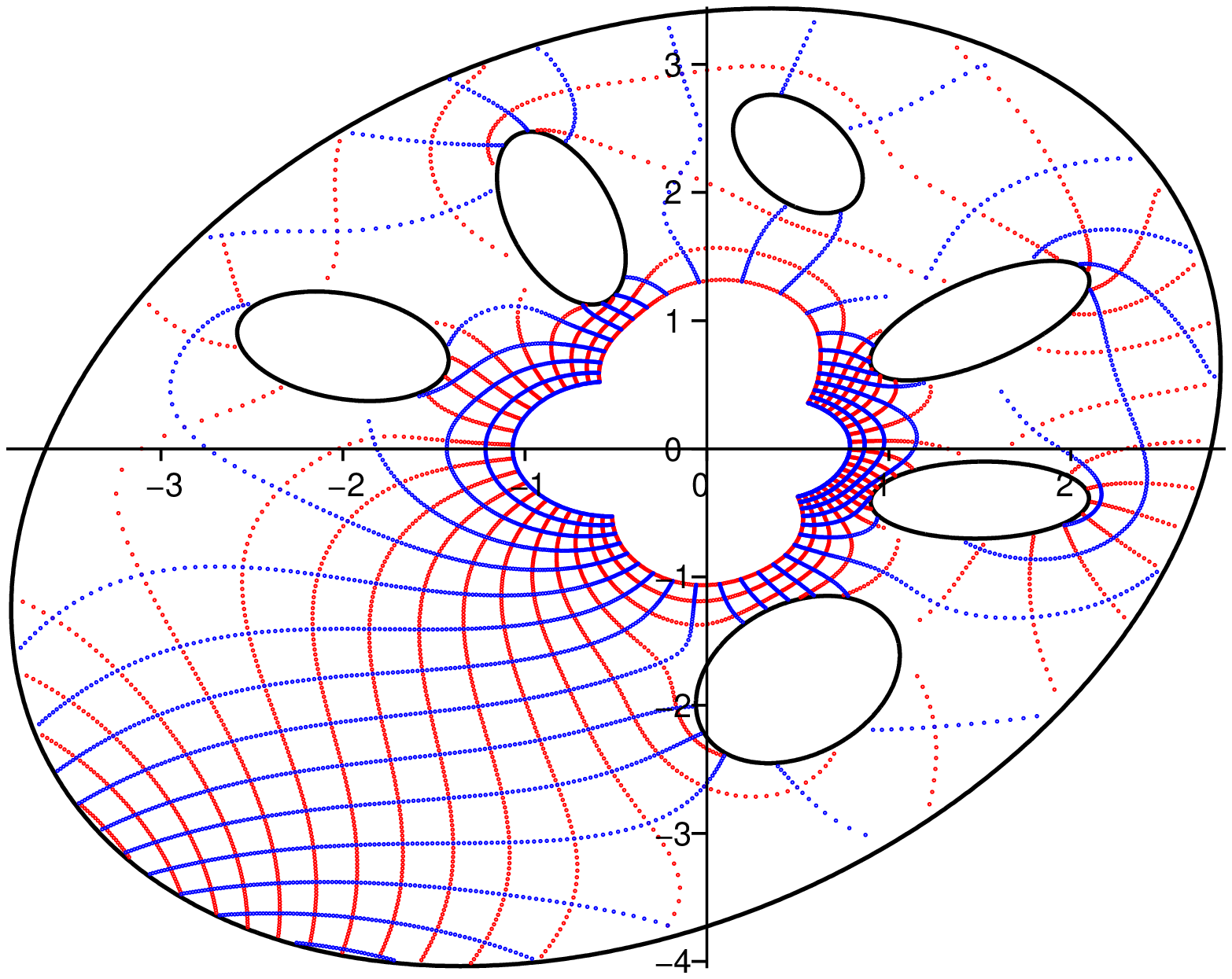}}%
\hfill}%
\caption{The inverse image of the circular slit canonical region.} 
\label{f:e2cr}
\end{figure}

\begin{figure}%
\centerline{%
\hfill\scalebox{0.25}{\includegraphics{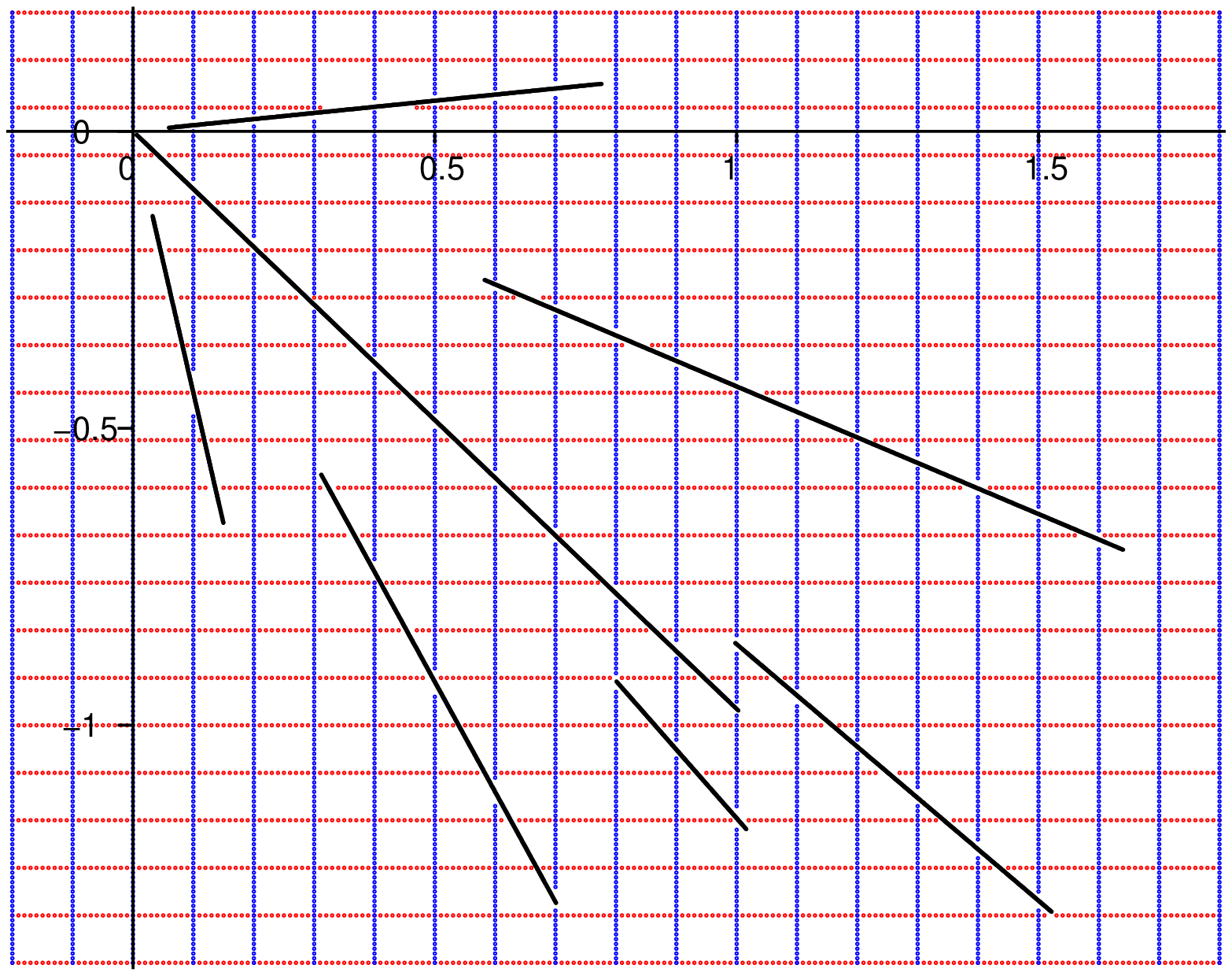}}%
\hfill\scalebox{0.25}{\includegraphics{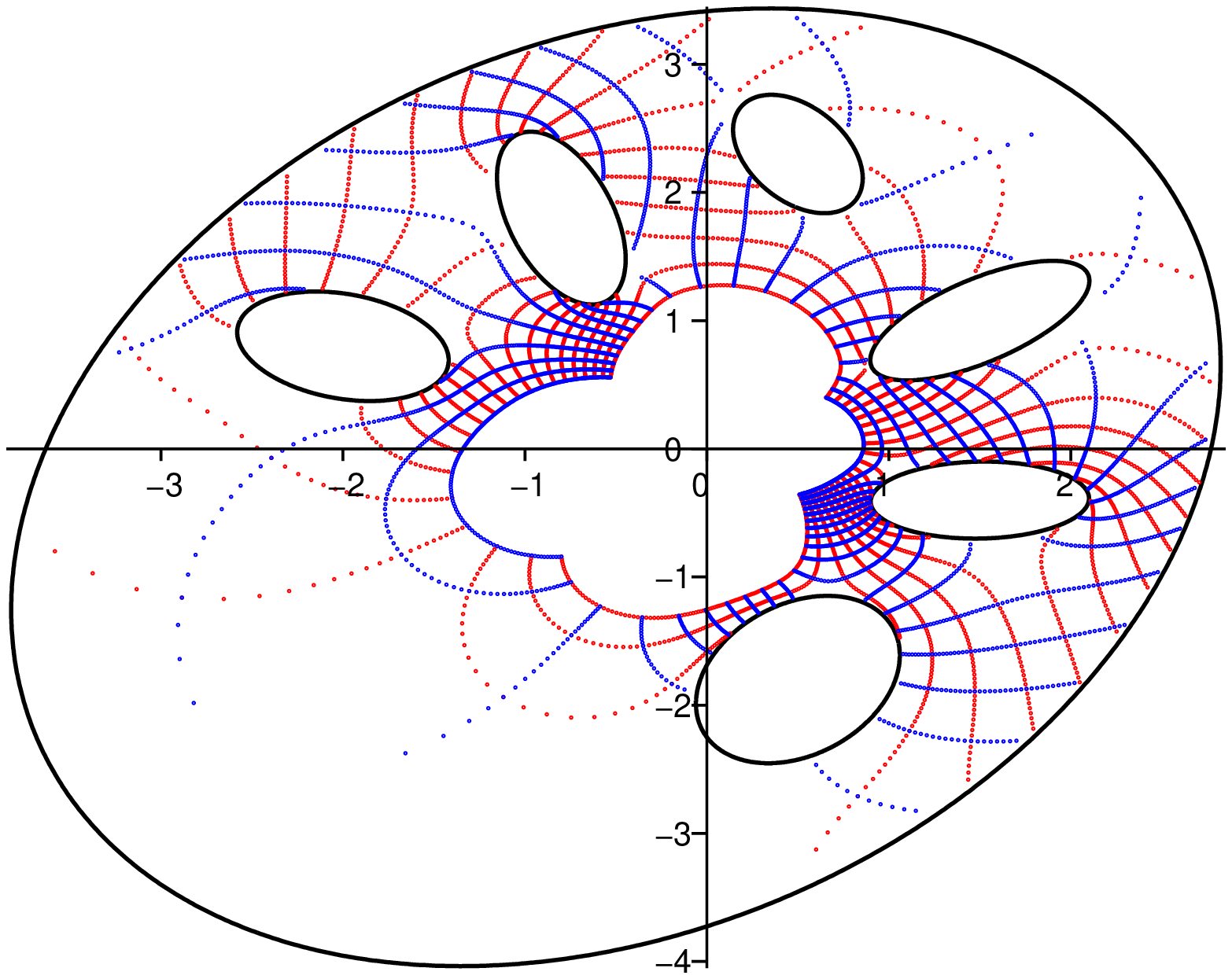}}%
\hfill}%
\caption{The inverse image of the radial slit canonical region.} 
\label{f:e2rd}
\end{figure}

\begin{figure}%
\centerline{%
\hfill\scalebox{0.25}{\includegraphics{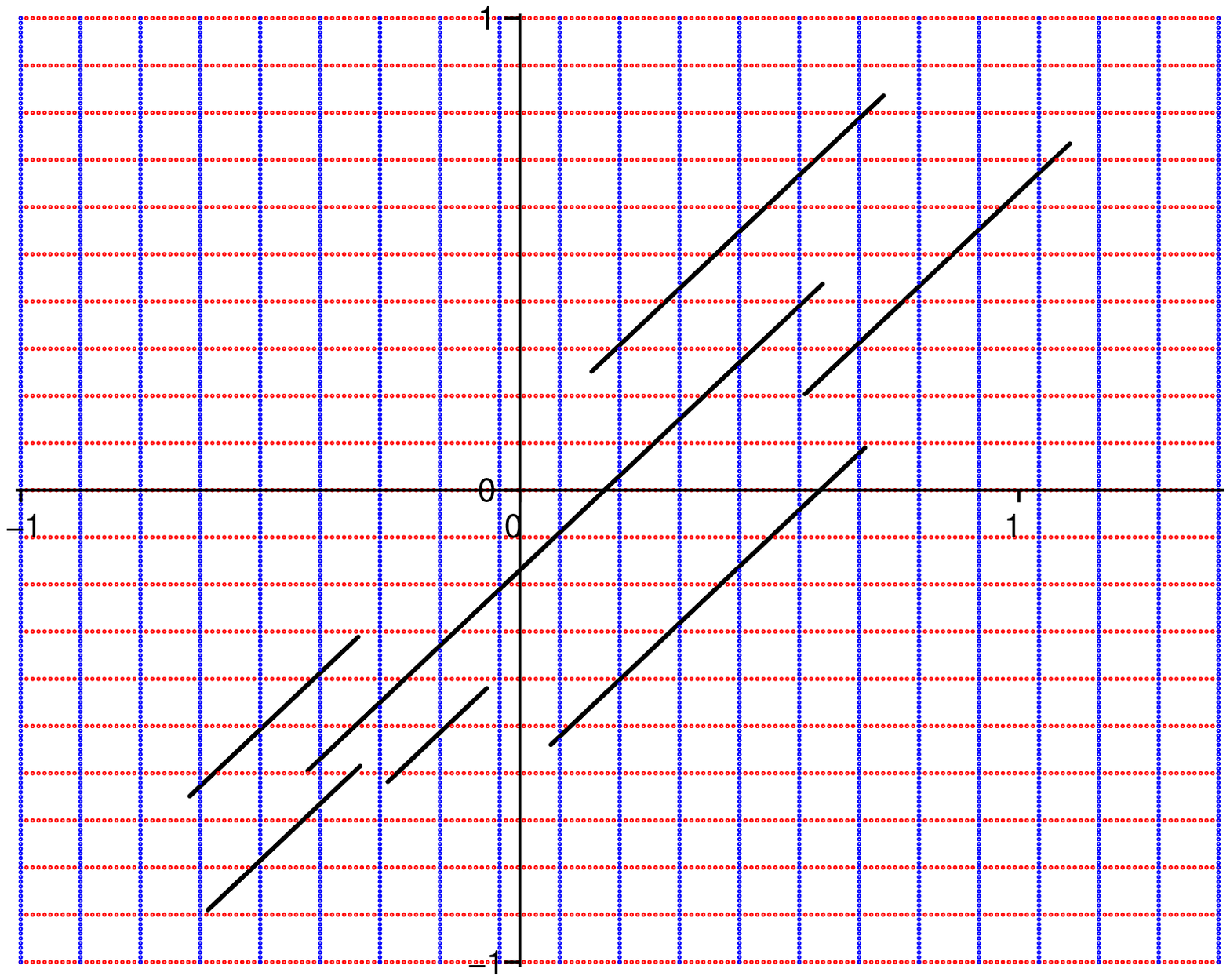}}%
\hfill\scalebox{0.25}{\includegraphics{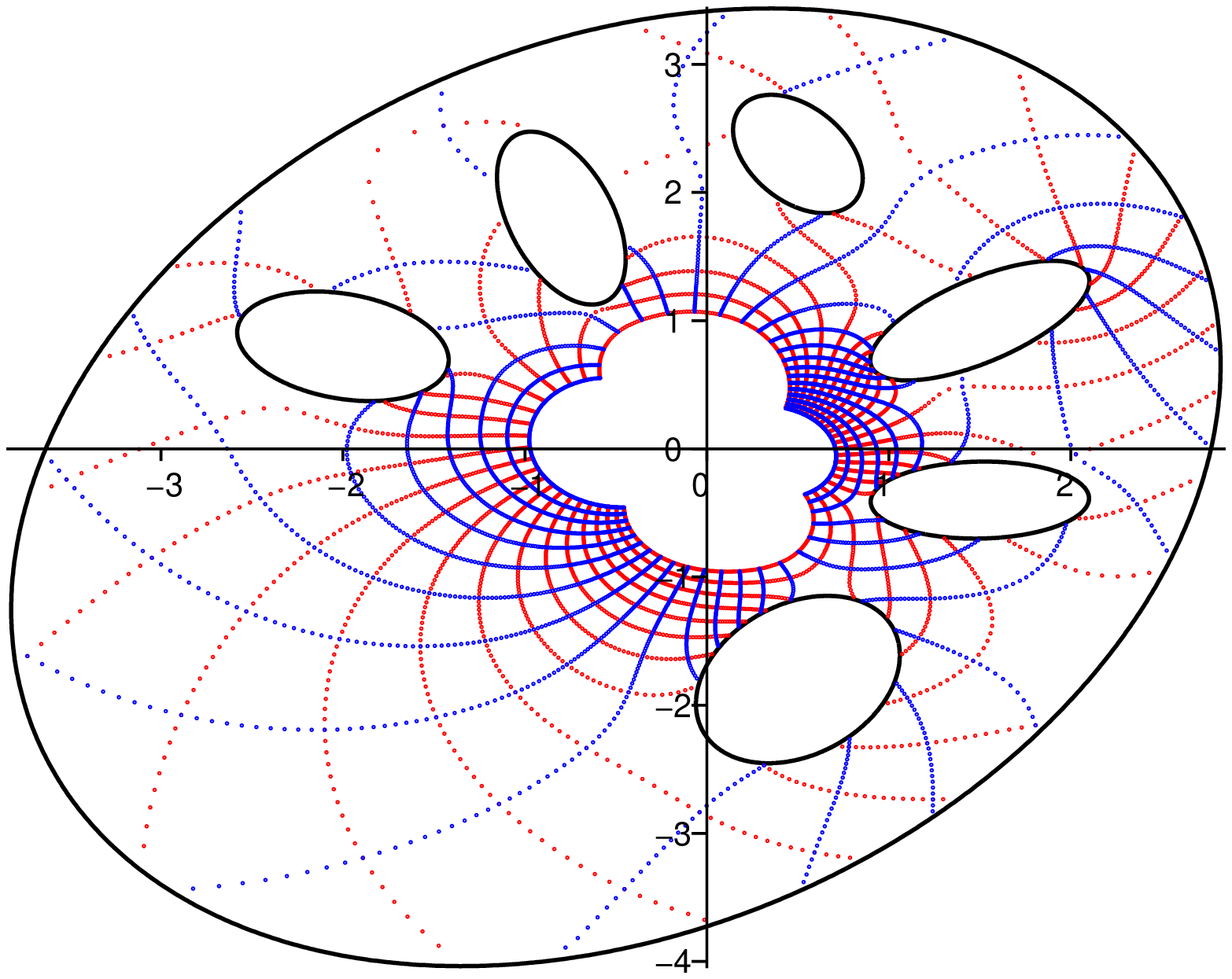}}%
\hfill}%
\caption{The inverse image of the parallel slit canonical region.} 
\label{f:e2pr}
\end{figure}

In this section, we consider a bounded multiply connected regions of connectivity $7$ (see Figure~\ref{f:e2or}(a)). The boundary $\Gamma$ of the bounded region $G$ is parametrized by
\begin{equation}\label{e:eta-ex}
\eta_j(t)=z_j+e^{\i\sigma_j}(\alpha_j\cos t+\i\beta_j\sin t), \quad
j=0,1,\ldots,6.
\end{equation}
The values of the complex constants $z_j$ and the real constants $\alpha_j$, $\beta_j$ and $\sigma_j$ are as in Table~\ref{t:cons}. We use the presented method to compute the mapping functions from the original region onto the five canonical slit regions and the inverse mapping functions from the five canonical regions onto the original region. For parallel silt region we set $\delta=\pi/4$. The numerical results, calculated with $n=1024$, are presented in Figures~\ref{f:e2or}--\ref{f:e2pr}. We present also a comparison between the presented method and the method presented in~\cite{nas-fun,nas-siam,nas-jmaa}. Table~\ref{t:error} lists the norm $\|\omega_n-\hat\omega_n\|_\infty$ for several values of $n$ where $\omega_n$ is the approximate boundary value of the mapping function obtained by the present method and $\hat\omega_n$ is the approximate boundary value of the mapping function obtained by the method presented in~\cite{nas-fun,nas-siam,nas-jmaa}.

\begin{table}
\caption{The values of constants $\alpha_j$, $\beta_j$, $z_j$ and $\sigma_j$ in~(\ref{e:eta-ex}).}
\label{t:cons}%
\vskip-0.25cm
\[
\begin{array}{r@{\hspace{0.75cm}}r@{\hspace{0.75cm}}r@{\hspace{0.75cm}}r@{\hspace{0.75cm}}r} \hline %
j   &  \alpha_j &  \beta_j & z_j\quad\quad & \sigma_j \\  %
\hline %

0  &4.0  & 3.0  &-0.5-0.3\i  &1.0 \\  
1  &0.7  &-0.3  & 1.5+1.0\i  &0.6 \\
2  &0.3  &-0.6  & 1.5-0.4\i  &1.6 \\
3  &0.5  &-0.7  & 0.5-1.8\i  &2.6 \\
4  &0.6  &-0.4  &-2.0+0.8\i  &2.8 \\
5  &0.3  &-0.7  &-0.8+1.8\i  &0.3 \\
6  &0.3  &-0.5  & 0.5+2.3\i  &0.5 \\
\hline %
\end{array}
\]
\end{table}

\begin{table}
\caption{The norm $\|\omega_n-\hat\omega_n\|_\infty$.}
\label{t:error}%
\vskip-0.25cm
\[
\begin{array}{l@{\hspace{0.25cm}}c@{\hspace{0.25cm}}c@{\hspace{0.25cm}}c@{\hspace{0.25cm}}c@{\hspace{0.25cm}}c} \hline %
n    &\text{An annulus} &\text{A disc with} &\text{The circular} &\text{The radial} &\text{The parallel} \\  %
     &\text{with circular slit} &\text{circular slit} &\text{slit region} &\text{slit region} &\text{slit region} \\  %
\hline %
16  &1.4(-01)&2.1(-01)&1.9(-01)&1.2(-01)&5.4(-02) \\  
32  &1.3(-02)&1.8(-02)&1.5(-02)&7.8(-03)&3.4(-03) \\
64  &1.6(-04)&1.7(-04)&1.6(-04)&1.3(-04)&5.2(-05) \\
128 &1.9(-07)&2.1(-07)&9.9(-08)&1.6(-07)&6.6(-08) \\
256 &5.0(-13)&6.2(-13)&2.8(-13)&3.8(-13)&1.8(-13) \\
\hline %
\end{array}
\]
\end{table}

\section{Conclusions}

A new uniquely solvable boundary integral equation has been presented in this paper for computing, in a unified way, the conformal mapping function $w=\omega(z)$, its derivative $\frac{dw}{dz}=\omega'(z)$, and its inverse $z=\omega^{-1}(w)$ from bounded multiply connected regions onto the five classical canonical slit regions. The presented boundary integral equation can also be used for simply connected regions. The integral equation has been derived by reformulating the conformal mapping as an adjoint RH problem. Then, based on the results presented in~\cite{wegm} and the properties of the boundary correspondence function $\theta(t)$, a uniquely solvable boundary integral equation has been derived for $\theta'(t)$. The function $\theta(t)$ was calculated by integrating the function $\theta'$ where the integration constants were computed using the same integral equation. The integral equation was used also to compute the parameters $R(t)=(R_0,\ldots,R_m)$ of the canonical regions. By obtaining $\theta(t)$ and $R(t)$, we obtain the boundary values of the mapping function $\omega(z)$ and its derivative $\omega'(z)$. The values of the mapping function $w=\omega(z)$ for $z\in G$ were calculated by means of the Cauchy integral formula. We computed also the values of the inverse mapping function $z=\omega^{-1}(w)$ by a Cauchy type integral that involves the functions $\theta(t)$, $\theta'(t)$ and $R(t)$. Since computing the functions $\theta(t)$, $\theta'(t)$ and $R(t)$ provides us with the boundary values of the mapping function $\omega(z)$ as well as the boundary values of its derivative $\omega'(z)$, we can also compute the values of the inverse mapping function by a Newton iteration method (see e.g.,~\cite[\S3.8]{kro}). 

The presented method is unified since it can be used to compute the mapping function, its derivative $\omega'(z)$, and its inverse onto the five classical canonical slit regions. Only the right-hand side of the integral equation is different from one canonical region to another. Computing the mapping function and its derivative at the same time is important in fluid dynamics (see e.g.,~\cite{bat,cro-uni,cro-lift,cro-slit,cro-kir,saf}).

The presented method has several advantages over the method presented in~\cite{nas-fun,nas-siam}. The right-hand side of the integral equation in~\cite{nas-fun,nas-siam} contains a singular operator which requires extra calculations. While, the right-hand side of the presented integral equation is given explicitly. The method presented in~\cite{nas-fun,nas-siam} was used to compute only the mapping function. However, the presented method can be used to compute the mapping function, its derivative, and its inverse.

The presented method also has several advantages over the method presented in~\cite{san-dc,san-pr,san-an,san-cr,san-rd}. Firstly, the presented method depends on only one integral equation which solvability is well known (see~\cite{nas-lap,wegm}). Whereas, the method presented in~\cite{san-dc,san-pr,san-an,san-cr,san-rd} depends on three different integral equations. The first of these three integral equation is different from a canonical region to another and its solvability has not been studied yet. Only numerical experiments have been given for this equation. Secondly, the present method is used to compute the mapping function as well as its inverse, while the method presented in~\cite{san-dc,san-pr,san-an,san-cr,san-rd} is used to compute only the mapping function. Finally, the boundary values of the mapping function $\omega(z)$ was computed from a formula that involves its derivative $\omega'(z)$, $\theta'$, and the parameters of the canonical regions. This formula contain the terms $\frac{\theta'(t)}{|\theta'(t)|}$. Since the function $\theta'(t)$ has positive and negative values, and it is zero at the end of the slits, it is difficult to calculate the values of the term $\frac{\theta'(t)}{|\theta'(t)|}$ in a stable way. When the values of $\theta'(t)$ are near zero, small changes in the values of $\theta'(t)$ could yield wrong values of the term $\frac{\theta'(t)}{|\theta'(t)|}$. We do not have such difficulties in the presented method.



\end{document}